\documentclass[11pt,reqno]{amsart}
\usepackage[latin1]{inputenc} 
\usepackage{indentfirst}
\usepackage{epsfig} 
\usepackage{graphicx} 
\usepackage{multirow} 
\usepackage{verbatim} 
\usepackage{rotating} 
\usepackage{lettrine}
\usepackage{lscape}
\usepackage{ae}
\usepackage{float}
\usepackage[all]{xy}
\usepackage{icomma} 
\usepackage{varioref} 
\usepackage{layout}
\usepackage[Lenny]{fncychap} 
\usepackage{enumerate}
\usepackage{amsfonts, amsmath, amssymb, stmaryrd, latexsym, amsthm, }
\usepackage{dsfont}
\usepackage{wrapfig, color, epsfig}
\usepackage{xspace, setspace}
\usepackage{array}
\usepackage{longtable}
\usepackage{geometry}
\usepackage[frenchb,english]{babel}
\begin{document}
\newtheorem{The}{Theorem}[section]

\newcommand{\K}{\mathbb{Q}(\sqrt{d},i)}
\newcommand{\q}{\mathbb{Q}(\sqrt{d})/\mathbb{Q}}
\newcommand{\qq}{\mathbb{Q}(\sqrt{d})}
\newcommand{\pro}{\prod_{p/\Delta_K}e(p)}
\newtheorem{lem}[The]{Lemma}
\newtheorem{propo}[The]{Proposition}
\newtheorem{coro}[The]{Corollary}
\newtheorem{proprs}[The]{properties}
\theoremstyle{remark}
\newtheorem{rema}[The]{\bf Remark}
\newtheorem{remas}[The]{\bf Remarks}
\newtheorem{exam}[The]{\textbf{Numerical Example}}
\newtheorem{exams}[The]{\textbf{Numerical Examples}}
\newtheorem{df}[The]{definition}
\newtheorem{dfs}[The]{definitions}
\def\NN{\mathds{N}}
\def\RR{\mathbb{R}}
\def\HH{I\!\! H}
\def\QQ{\mathbb{Q}}
\def\CC{\mathbb{C}}
\def\ZZ{\mathbb{Z}}
\def\OO{\mathcal{O}}
\def\kk{\mathds{k}}
\def\KK{\mathbb{K}}
\def\FF{\mathbb{F}}
\def\ho{\mathcal{H}_0^{\frac{h(d)}{2}}}
\def\LL{\mathbb{L}}
\def\o{\mathcal{H}_0}
\def\h{\mathcal{H}_1}
\def\hh{\mathcal{H}_2}
\def\hhh{\mathcal{H}_3}
\def\hhhh{\mathcal{H}_4}
\def\k{\mathds{k}^{(*)}}
\def\G{\mathds{k}^{(*)}}
\def\l{\mathds{L}}
\def\L{\kk_2^{(2)}}
\def\kkk{k^{(*)}}
\def\Q{\mathcal{Q}}
\def\D{\kk_2^{(1)}}
\title[ On the strongly ambiguous...]{On the strongly ambiguous classes of some biquadratic number fields}
\author[Abdelmalek AZIZI]{Abdelmalek Azizi}
\address{Abdelmalek Azizi and Abdelkader Zekhnini: Département de Mathématiques, Faculté des Sciences, Université Mohammed 1, Oujda, Morocco }
\author{Abdelkader Zekhnini}
\email{abdelmalekazizi@yahoo.fr}
\email{zekha1@yahoo.fr}
\author[Mohammed Taous]{Mohammed Taous}
\address{Mohammed Taous: Département de Mathématiques, Faculté des Sciences et Techniques, Université Moulay Ismail, Errachidia, Morocco}
\email{taousm@hotmail.com}

\maketitle
\selectlanguage{english}
\begin{abstract}
We study the capitulation of ideal classes in an infinite
family of imaginary bicyclic biquadratic number fields consisting of
 fields $\kk =\QQ(\sqrt{2pq}, i)$, where $i=\sqrt{-1}$ and $p\equiv -q\equiv1 \pmod 4$
are different primes. For each of the three quadratic extensions $\KK/\kk$
inside the absolute genus field $\kk^{(*)}$ of $\kk$, we compute the capitulation
kernel of $\KK/\kk$. Then we deduce that each strongly ambiguous class of $\kk/\QQ(i)$
capitulates already in $\kk^{(*)}$, which is  smaller than the relative
genus field $(\kk/\QQ(i))^*$.
\end{abstract}
\section{Introduction}\label{sec1}
Let $k$ be an algebraic number field and let $\mathbf{C}l_2(k)$ denote its 2-class group, that is the 2-Sylow
subgroup of the ideal class group, $\mathbf{C}l(k)$,  of $k$.
We denote by  $k^{(*)}$  the absolute genus field of $k$.
Suppose $F$ is a finite extension of $k$,
then we say that an ideal class of $k$ capitulates in $F$ if it is in the kernel of the
homomorphism
$$J_F: \mathbf{C}l(k) \longrightarrow \mathbf{C}l(F)$$ induced by extension of ideals from $k$
to $F$. An
important problem in Number Theory is to determine explicitly the kernel of  $J_F$, which is
usually called the capitulation kernel. If $F$ is
the  relative genus  field of a cyclic extension $K/k$, which we denote  by $(K/k)^*$ and that is  the maximal unramified extension  of  $K$ which is obtained by composing  $K$ and an abelian extension over $k$, F. Terada  states in  \cite{FT-71} that all the ambiguous ideal classes of $K/k$, which are classes of $K$ fixed under any element of $\mathrm{G}al(K/k)$,  capitulate in
 $(K/k)^*$. If $F$ is the  absolute genus  field of an abelian extension $K/\QQ$, then H. Furuya confirms in \cite{Fu-77} that every   strongly ambiguous  class of $K/\QQ$, that is an ambiguous ideal class containing at least one ideal invariant under any element of $\mathrm{G}al(K/\QQ)$, capitulates in  $F$. In this paper, we construct a family of number fields $k$ for which $\mathbf{C}l_2(k)\simeq (2, 2, 2)$ and all the strongly ambiguous  classes of $k/\QQ(i)$ capitulate in  $k^{(*)}\varsubsetneq (k/\QQ(i))^*$.\par
Let $p$ and $q$ be different primes,  $\mathds{k}=\QQ(\sqrt{2pq},i)$ and $\mathds{K}$ be an unramified quadratic extension of $\kk$ that is abelian  over $\QQ$. Denote by $\mathrm{A}m_s(\kk/\QQ(i))$
 the group of the strongly ambiguous classes of  $\kk/\QQ(i)$. In \cite{Az-00}, the  first author has studied the capitulation problem in $\KK/\kk$ assuming $p\equiv -q\equiv1\pmod4$ and  $\mathbf{C}l_2(\kk)\simeq (2, 2)$. On the other hand, in \cite{AZT14-1},  we have dealt with  the same problem assuming $p\equiv q\equiv1\pmod4$,  and in \cite{AZT14-2}, we have studied the capitulation problem of the $2$-ideal classes of $\kk$ in its fourteen unramified extensions, within the first Hilbert $2$-class field of $\kk$,  assuming $p\equiv q\equiv5\pmod8$. It is the purpose of the present article
to pursue this research project further  for all types of $\mathbf{C}l_2(\kk)$, assuming $p\equiv -q\equiv1\pmod4$,  we compute the capitulation
kernel of $\KK/\kk$ and  we deduce  that   $\mathrm{A}m_s(\kk/\QQ(i))\subseteq \ker J_{\k}$. As
an application we will determine these kernels when  $\mathbf{C}l_2(\mathds{k})$ is of type $(2, 2, 2)$.

Let  $k$ be a number field, during this paper, we adopt the following notations:
 \begin{itemize}
   \item   $p\equiv -q\equiv1\pmod4$ are different primes.
   \item $\kk$: denotes the field $\QQ(\sqrt{2pq},\sqrt{-1})$.
   \item $\kappa_{K}$: the capitulation kernel of an unramified extension $K/\kk$.
   \item $\mathcal{O}_k$: the ring of integers of $k$.
   \item $E_k$: the unit group of $\mathcal{O}_k$.
   \item  $W_k$: the group of roots of unity contained in $k$.
   \item $\mathrm{F}.\mathrm{S}.\mathrm{U}$ : the fundamental system of units.
   \item $k^+$: the maximal real subfield of $k$, if it is a CM-field.
   \item $Q_k=[E_k:W_kE_{k^+}]$ is Hasse's unit index, if $k$ is a CM-field.
   \item $q(k/\QQ)=[E_k:\prod_{i=1}^{s} E_{k_i}]$ is the unit index of $k$, if $k$ is multiquadratic, where $k_1$, ..., $k_s$  are the  quadratic subfields of $k$.
  \item $k^{(*)}$: the absolute genus field of $k$.
  \item $\mathbf{C}l_2(k)$: the 2-class group of $k$.
  \item $i=\sqrt{-1}$.
  \item $\epsilon_m$: the fundamental unit of $\QQ(\sqrt m)$, if $m>1$ is a square-free integer.
  \item $N(a)$: denotes the absolute norm of a number $a$ i.e. $N_{k/\QQ}(a)$, where $k=\QQ(\sqrt a)$.
  \item $x\pm y$ means $x+y$ or $x-y$ for some numbers $x$ and $y$.
 \end{itemize}
\section{\bf{Preliminary results}}
Let us  first collect some results that will be useful in what follows.\par
Let $k_j$, $1\leq j\leq 3$, be the three real quadratic subfields of a biquadratic bicyclic real number field $K_0$
and $\epsilon_j>1$ be the fundamental unit of $k_j$. Since $\alpha^2N_{K_0/{\bf Q}}(\alpha) =\prod_{j=1}^3N_{K_0/k_j}(\alpha)$ for any $\alpha\in K_0$,
the square of any unit of $K_0$ is in the group generated by the $\epsilon_j$'s, $ 1\leq j\leq 3$.
Hence, to determine a system of fundamental units of $K_0$ it suffices to determine which of the units in
$B:=\{\epsilon_1,\epsilon_2,\epsilon_3,\epsilon_1\epsilon_2,\epsilon_1\epsilon_3,\epsilon_2\epsilon_3,\epsilon_1\epsilon_2\epsilon_3\}$
are  squares in $K_0$ (see \cite{Wa-66} or \cite{Lo-00}).
Put  $K=K_0(i)$,  then  to determine a $\mathrm{F}.\mathrm{S}.\mathrm{U}$  of $K$,  we will use the following result (see \cite[p.18]{Az-99})  that the first author has deduced from a theorem of Hasse \cite[\S 21, Satz 15 ]{Ha-52}.
\begin{lem}\label{3}
Let $n\geq2$ be an integer   and $\xi_n$ a $2^n$-th primitive root of unity, then
$$\begin{array}{lllr}\xi_n = \dfrac{1}{2}(\mu_n + \lambda_ni), &\hbox{ where }&\mu_n =\sqrt{2 +\mu_{n-1}},&\lambda_n =\sqrt{2 -\mu_{n-1}}, \\
                               &    & \mu_2=0, \lambda_2=2 &\hbox{ and\quad } \mu_3=\lambda_3=\sqrt 2.
                                              \end{array}
                                            $$
Let $n_0$ be the greatest integer such that  $\xi_{n_0}$ is contained in $K$, $\{\epsilon'_1, \epsilon'_2, \epsilon'_3\}$ a $\mathrm{F}.\mathrm{S}.\mathrm{U}$  of $K_0$ and $\epsilon$ a unit of  $K_0$ such that $(2 + \mu_{n_0})\epsilon$ is a square in  $K_0$ $($if it exists$)$. Then a $\mathrm{F}.\mathrm{S}.\mathrm{U}$  of  $K$ is one of the following systems:
\begin{enumerate}[\rm(1)]
\item $\{\epsilon'_1, \epsilon'_2, \epsilon'_3\}$ if $\epsilon$ does not exist,
\item $\{\epsilon'_1, \epsilon'_2, \sqrt{\xi_{n_0}\epsilon}\}$ if $\epsilon$
exists;  in this case $\epsilon = {\epsilon'_1}^{i_1}{
\epsilon'_2}^{i_2}\epsilon'_3$, where $i_1$, $i_2\in \{0, 1\}$ (up to a permutation).
\end{enumerate}
\end{lem}
\begin{lem}[\cite{Az-00}, {Lemma 5}]\label{5}
Let $d>1$ be a square-free integer and $\epsilon_d=x+y\sqrt d$,
where $x$, $y$ are  integers or semi-integers. If $N(\epsilon_d)=1$, then $2(x+1)$, $2(x-1)$, $2d(x+1)$ and
$2d(x-1)$ are not squares in  $\QQ$.
\end{lem}
\begin{lem}[\cite{Az-00}, {Lemma 6}]\label{1:048}
Let $q\equiv-1\pmod4$ be a prime and $\epsilon_q=x+y\sqrt q$ be the fundamental unit of  $\QQ(\sqrt q)$. Then $x$ is an even integer, $x\pm1$ is a square in  $\NN$ and $2\epsilon_q$ is a square in $\QQ(\sqrt q)$.
\end{lem}
\begin{lem}[\cite{Az-00}, {Lemma 7}] \label{4}
Let  $p$ be an odd  prime and $\epsilon_{2p}=x+y\sqrt{2p}$.
If  $N(\epsilon_{2p})=1$, then $x\pm1$ is a square in $\NN$ and $2\epsilon_{2p}$
is a square in $\QQ(\sqrt{2p})$.
\end{lem}
\begin{lem}[\cite{Az-99}, {3.(1) p.19}]\label{6}
Let $d>2$ be a square-free integer and $k=\QQ(\sqrt d,i)$, put $\epsilon_d=x+y\sqrt d$.
\begin{enumerate}[\rm\indent(1)]
  \item If $N(\epsilon_d)=-1$, then $\{\epsilon_d\}$ is a $\mathrm{F}.\mathrm{S}.\mathrm{U}$  of $k$.
  \item If $N(\epsilon_d)=1$, then $\{\sqrt{i\epsilon_d}\}$ is a $\mathrm{F}.\mathrm{S}.\mathrm{U}$  of $k$ if
   and only if $x\pm1$ is a square in $\NN$ i.e. $2\epsilon_d$ is a square in $\QQ(\sqrt d)$. Else $\{\epsilon_d\}$ is a $\mathrm{F}.\mathrm{S}.\mathrm{U}$  of $k$ $($this result is also in \cite{Kub-56}$)$.
\end{enumerate}
\end{lem}
\section{\textbf{$\mathrm{F}.\mathrm{S}.\mathrm{U}$  OF SOME CM-FIELDS}}
As  $\kk=\QQ(\sqrt{2pq}, i)$, so $\kk$ admits three unramified quadratic extensions that are abelian over $\QQ$, which are $\KK_1=\kk(\sqrt{p})=\QQ(\sqrt{p}, \sqrt{2q}, i)$,   $\KK_2=\kk(\sqrt{q})=\QQ(\sqrt{q}, \sqrt{2p}, i)$ and $\KK_3=\kk(\sqrt{2})=\QQ(\sqrt{2}, \sqrt{pq}, i)$. Put $\epsilon_{2pq}=x+y\sqrt{2pq}$, the  first author has given in \cite{Az-00} the $\mathrm{F}.\mathrm{S}.\mathrm{U}$'s of these three fields, if $2\epsilon_{2pq}$  is not a square in $\QQ(\sqrt{2pq})$ i.e. $x+1$ and $x-1$  are not squares in $\NN$. In what follows, we determine the $\mathrm{F}.\mathrm{S}.\mathrm{U}$'s of $\KK_j$, $1\leq j\leq3$, in all cases.
\subsection{\textbf{$\mathrm{F}.\mathrm{S}.\mathrm{U}$  of the field $\KK_1$}}
Let $\KK_1=\kk(\sqrt{p})=\QQ(\sqrt{p}, \sqrt{2q}, i)$.
\begin{propo}\label{27}
Keep the previous notations, then $Q_{\KK_1}=2$ and just one of the following two cases holds:
\begin{enumerate}[\rm\indent(1)]
  \item If $x\pm1$ or $p(x\pm1)$ is a square in $\NN$, then $\left\{\epsilon_p, \epsilon_{2q}, \sqrt{\epsilon_{2q}\epsilon_{2pq}}\right\}$ is a  $\mathrm{F}.\mathrm{S}.\mathrm{U}$ of $\KK_1^+$ and that of $\KK_1$ is $\left\{\epsilon_{p}, \sqrt{i\epsilon_{2q}}, \sqrt{\epsilon_{2q}\epsilon_{2pq}}\right\}$.
  \item If $2p(x\pm1)$ is a square in $\NN$, then $\left\{\epsilon_p, \epsilon_{2q}, \sqrt{\epsilon_{2pq}}\right\}$ is a $\mathrm{F}.\mathrm{S}.\mathrm{U}$ of  $\KK_1^+$ and that of $\KK_1$ is  $\left\{\epsilon_p, \sqrt{i\epsilon_{2q}}, \sqrt{\epsilon_{2pq}}\right\}$.
\end{enumerate}
\end{propo}
\begin{proof}
As  $p\equiv1\pmod4$, then $\epsilon_p$ is not a square in $\KK_1^+$; but   $\epsilon_{2pq}$ and $\epsilon_{2q}\epsilon_{2pq}$ can be. Moreover, according to Lemma \ref{4},  $2\epsilon_{2q}$ is a square in  $\QQ(\sqrt{2q})$. On the other hand, we know that $N(\epsilon_{2pq})=1$, then $(x\pm1)(x\mp1)=2pqy^2$. Hence, by Lemma \ref{5} and according to the decomposition uniqueness in $\ZZ$, there are  three possibilities: $x\pm1$ or $p(x\pm1)$ or $2p(x\pm1)$ is a square in $\NN$, the only remaining case is the first one.
 If $x\pm1$ is a square in $\NN$ (for the other cases see \cite{Az-00}), then, by Lemma \ref{6}, $2\epsilon_{2pq}$ is a square in $\KK_1$.
Consequently,   $\sqrt{\epsilon_{2q}\epsilon_{2pq}}\in\KK_1^+$; hence  $\left\{\epsilon_p, \epsilon_{2q}, \sqrt{\epsilon_{2q}\epsilon_{2pq}}\right\}$ is a  $\mathrm{F}.\mathrm{S}.\mathrm{U}$ of $\KK_1^+$, and since  $2\epsilon_{2q}$ is a square in $\KK_1^+$, so Lemma \ref{3} yields that $\left\{\epsilon_p, \sqrt{i\epsilon_{2q}}, \sqrt{\epsilon_{2q}\epsilon_{2pq}}\right\}$ is a  $\mathrm{F}.\mathrm{S}.\mathrm{U}$ of $\KK_1$. Thus $Q_{\KK_1}=2$.
\end{proof}
\subsection{\textbf{$\mathrm{F}.\mathrm{S}.\mathrm{U}$  of the field  $\KK_2$}}
Let $\KK_2=\kk(\sqrt{q})=\QQ(\sqrt{q}, \sqrt{2p}, i)$.
\begin{propo}\label{28}
Keep the previous notations, then $Q_{\KK_2}=2$.
\begin{enumerate}[\rm\indent(1)]
\item Assume that $N(\epsilon_{2p})=1$. Then  just one of the following two cases holds.
\begin{enumerate}[\rm(i)]
  \item If $x\pm1$ or $2p(x\pm1)$ is a square in $\NN$, then $\left\{\sqrt{\epsilon_q\epsilon_{2p}}, \sqrt{\epsilon_q\epsilon_{2pq}}, \sqrt{\epsilon_{2p}\epsilon_{2pq}}\right\}$ is a  $\mathrm{F}.\mathrm{S}.\mathrm{U}$ of $\KK_2^+$ and that of $\KK_2$ is $\left\{\sqrt{i\epsilon_q}, \sqrt{i\epsilon_{2p}}, \sqrt{i\epsilon_{2pq}}\right\}$.
  \item If  $p(x\pm1)$ is a square in $\NN$, then $\left\{\epsilon_q, \sqrt{\epsilon_q\epsilon_{2p}}, \sqrt{\epsilon_{2pq}}\right\}$ is a  $\mathrm{F}.\mathrm{S}.\mathrm{U}$ of  $\KK_2^+$ and that of $\KK_2$ is  $\left\{\sqrt{i\epsilon_q}, \sqrt{i\epsilon_{2p}}, \sqrt{\epsilon_{2pq}}\right\}$.
\end{enumerate}
\item Assume that $N(\epsilon_{2p})=-1$. Then  just one of the following two cases holds.
\begin{enumerate}[\rm(i)]
  \item If $x\pm1$ or $2p(x\pm1)$ is a square in $\NN$, then $\left\{\epsilon_q, \epsilon_{2p}, \sqrt{\epsilon_q\epsilon_{2pq}}\right\}$ is a  $\mathrm{F}.\mathrm{S}.\mathrm{U}$ of $\KK_2^+$ and that of $\KK_2$ is  $\left\{\sqrt{i\epsilon_q}, \epsilon_{2p}, \sqrt{\epsilon_q\epsilon_{2pq}}\right\}$.
  \item If   $p(x\pm1)$ is a square in $\NN$, then $\left\{\epsilon_q, \epsilon_{2p},  \sqrt{\epsilon_{2pq}}\right\}$ is a  $\mathrm{F}.\mathrm{S}.\mathrm{U}$ of $\KK_2^+$ and that of $\KK_2$ is  $\left\{\sqrt{i\epsilon_q}, \epsilon_{2p}, \sqrt{\epsilon_{2pq}}\right\}$.
\end{enumerate}
\end{enumerate}
\end{propo}
\begin{proof}
According to Lemma \ref{6}, if  $x\pm1$ is a square in $\NN$, then  $2\epsilon_{2pq}$ is a square in $\QQ(\sqrt{2pq})$. Moreover, Lemma \ref{1:048} implies that $2\epsilon_{q}$ is also a square in $\QQ(\sqrt{q})$.
\begin{enumerate}
\item If $N(\epsilon_{2p})=1$, then  Lemma \ref{4} yields that $2\epsilon_{2p}$ is a square in $\QQ(\sqrt{2p})$, thus $\epsilon_{2p}\epsilon_{2pq}$, $\epsilon_{q}\epsilon_{2pq}$ and $\epsilon_{q}\epsilon_{2p}$ are squares in  $\KK_2^+$,  this gives the $\mathrm{F}.\mathrm{S}.\mathrm{U}$ of $\KK_2^+$, and that of $\KK_2$ is deduced by Lemma  \ref{3}.
\item If $N(\epsilon_{2p})=-1$, then   $\epsilon_{q}\epsilon_{2pq}$ is a square in  $\KK_2^+$,  this gives the $\mathrm{F}.\mathrm{S}.\mathrm{U}$ of $\KK_2^+$, and that of $\KK_2$ is deduced by Lemma  \ref{3}.
\end{enumerate}
For the other cases see \cite{Az-00}.
\end{proof}
\subsection{\textbf{$\mathrm{F}.\mathrm{S}.\mathrm{U}$  of the field  $\KK_3$}}
Let $\KK_3=\kk(\sqrt{2})=\QQ(\sqrt{2}, \sqrt{pq}, i)$.
\begin{propo}\label{29}
Put $\epsilon_{pq}=a+b\sqrt{pq}$, where $a$ and $b$ are in $\ZZ$.
\begin{enumerate}[\rm\indent(1)]
  \item If both of $x\pm1$ and  $a\pm1$ are squares in  $\NN$, then
\begin{enumerate}[\rm(i)]
 \item If $Q_{\KK_3}=1$, then $\left\{\epsilon_{2}, \sqrt{\epsilon_{pq}}, \sqrt{\epsilon_{2pq}}\right\}$ is a  $\mathrm{F}.\mathrm{S}.\mathrm{U}$ of both $\KK_3^+$ and $\KK_3$.
\item If $Q_{\KK_3}=2$, then $\left\{\epsilon_{2}, \sqrt{\epsilon_{pq}}, \sqrt{\epsilon_{2pq}}\right\}$ is a  $\mathrm{F}.\mathrm{S}.\mathrm{U}$ of  $\KK_3^+$ and that of $\KK_3$ is  $\left\{\epsilon_{2}, \sqrt{\epsilon_{pq}}, \sqrt{\xi\sqrt{\epsilon_{pq}\epsilon_{2pq}}}\right\}$, where $\xi$ is  an 8-th root of unity.
\end{enumerate}
  \item If  $x\pm1$ is a square in $\NN$ and  $a+1$, $a-1$ are not, then $\left\{\epsilon_{2}, \epsilon_{pq}, \sqrt{\epsilon_{2pq}}\right\}$ is a  $\mathrm{F}.\mathrm{S}.\mathrm{U}$ of both $\KK_3^+$ and $\KK_3$; hence $Q_{\KK_3}=1$.
  \item If  $a\pm1$ is a square in $\NN$ and  $x+1$, $x-1$ are not, then $\left\{\epsilon_{2}, \epsilon_{2pq},   \sqrt{\epsilon_{pq}}\right\}$ is a  $\mathrm{F}.\mathrm{S}.\mathrm{U}$ of both $\KK_3^+$ and $\KK_3$; hence $Q_{\KK_3}=1$.
  \item If   $x+1$, $x-1$, $a+1$ and $a-1$ are not squares in $\NN$, then $\left\{\epsilon_{2}, \epsilon_{pq}, \sqrt{\epsilon_{pq}\epsilon_{2pq}}\right\}$ is a $\mathrm{F}.\mathrm{S}.\mathrm{U}$ of both $\KK_3^+$ and $\KK_3$; hence $Q_{\KK_3}=1$.
\end{enumerate}
\end{propo}
Before proving this proposition, we quote the following  result.
\begin{rema}\label{37}
Keep the notations and hypotheses  of the Proposition $\ref{29}$.
\begin{enumerate}[\rm(1)]
 \item  if at most one of the numbers  $x+1$,  $x-1$,  $a+1$ and  $a-1$ is  a square in  $\NN$, then according to \cite[Remarque $13$, p. $391$ ]{Az-00},  $\KK_3^+$ and  $\KK_3$ have the same  $\mathrm{F}.\mathrm{S}.\mathrm{U}$.
\item From \cite[Theorem $2$, p. $348$ ]{HiYo-90}, if both of $x\pm1$ and  $a\pm1$ are squares in  $\NN$, then the unit index of  $\KK_3$ is  $1$ or $2$.
\end{enumerate}
\end{rema}
\begin{proof}
We know that  $N(\epsilon_{2})=-1$  and $N(\epsilon_{pq})=N(\epsilon_{2pq})=1$.  Moreover,  $(2+\sqrt2)\epsilon_{2}^i\epsilon_{pq}^j\epsilon_{2pq}^k$ cannot be a square in $\KK_3^+$, for all $i$, $j$ and $k$ of $\left\{0, 1\right\}$; as otherwise with some   $\alpha\in\KK_3^+$ we would have $\alpha^2=(2+\sqrt 2)\epsilon_{2}^i\epsilon_{pq}^j\epsilon_{2pq}^k$, so $(N_{\KK_3^+/\QQ(\sqrt{pq})}(\alpha))^2=2(-1)^{i}\epsilon_{pq}^{2j}$, yielding  that  $\sqrt{\pm2}\in\QQ(\sqrt{pq})$, which is absurd.\\
\indent   As $a^2-1=pqb^2$, so  by Lemma \ref{5} and according to the decomposition uniqueness in $\ZZ$, there are  three possible cases: $a\pm1$ or $p(a\pm1)$ or $2p(a\pm1)$ is a square in $\NN$.
\begin{enumerate}[(a)]
  \item If $a\pm1$ is a square in $\NN$, then there exist $b_1$ and $b_2$ in $\NN$ with $b=b_1b_2$ such that
 $\left\{\begin{array}{ll}
  a\pm1=b_1^2,\\
  a\mp1=pqb_2^2,
  \end{array}\right.\text{ hence }\sqrt{\epsilon_{pq}}=\frac{1}{2}(b_1\sqrt2+b_2\sqrt{2pq})\in\KK_3^+.$
  \item If $p(a\pm1)$  is a square in $\NN$, then there exist $b_1$ and $b_2$ in $\NN$ with $b=b_1b_2$ such that
 $\left\{\begin{array}{ll}
  a\pm1=pb_1^2,\\
  a\mp1=qb_2^2,
  \end{array}\right. \text{ hence }
  \left\{\begin{array}{ll}
  \sqrt{\epsilon_{pq}}=\frac{1}{2}(b_1\sqrt{2p}+b_2\sqrt{2q})\not\in\KK_3^+,\\
  \sqrt{p\epsilon_{pq}}\in\KK_3^+ \ and \  \sqrt{q\epsilon_{pq}}\in\KK_3^+.
  \end{array}\right.$
  \item If $2p(a\pm1)$  is a square in $\NN$, then there exist $b_1$ and $b_2$ in $\NN$ with $b=b_1b_2$ such that
  $\left\{\begin{array}{ll}
  a\pm1=2pb_1^2,\\
  a\mp1=2qb_2^2
  \end{array}\right. \text{ hence }
  \left\{\begin{array}{ll}
  \sqrt{\epsilon_{pq}}=b_1\sqrt{p}+b_2\sqrt{q}\not\in\KK_3^+;\\
  \sqrt{p\epsilon_{pq}}\in \KK_3^+\  and \ \sqrt{q\epsilon_{pq}}\in \KK_3^+.
  \end{array}\right.$
\end{enumerate}
 Similarly, we get:
 \begin{enumerate}[(a')]
   \item If $x\pm1$ is a square in $\NN$,  then  $\sqrt{\epsilon_{2pq}}\in\KK_3^+$.
   \item If $p(x\pm1)$ is a square in $\NN$, then $\sqrt{\epsilon_{2pq}}\not\in\KK_3^+$, $\sqrt{p\epsilon_{2pq}}\in\KK_3^+$  and $\sqrt{q\epsilon_{2pq}}\in\KK_2^+$.
   \item If $2p(x\pm1)$ is a square in $\NN$, then $\sqrt{\epsilon_{2pq}}\not\in\KK_3^+$,  $\sqrt{p\epsilon_{2pq}}\in \KK_2^+$  and  $\sqrt{q\epsilon_{2pq}}\in \KK_2^+$.
 \end{enumerate}
Consequently, we find:
\begin{enumerate}[\rm(1)]
  \item If $a\pm1$ and $x\pm1$ are squares in $\NN$, then $\left\{\epsilon_{2}, \sqrt{\epsilon_{pq}}, \sqrt{\epsilon_{2pq}}\right\}$ is a $\mathrm{F}.\mathrm{S}.\mathrm{U}$ of $\KK_3^+$.
\begin{enumerate}[\rm(i)]
   \item If $Q_{\KK_3}=1$, then  $\left\{\epsilon_{2}, \sqrt{\epsilon_{pq}}, \sqrt{\epsilon_{2pq}}\right\}$ is also a  $\mathrm{F}.\mathrm{S}.\mathrm{U}$   of  $\KK_3$.
   \item If $Q_{\KK_3}=2$, then, according to \cite{HiYo-90}, $\KK_3^+(\sqrt{2+\sqrt2})=\KK_3^+(\sqrt{\sqrt{\epsilon_{pq}\epsilon_{2pq}}})$, so there exist  $\alpha\in\KK_3^+$ such that $2+\sqrt2=\alpha^2\sqrt{\epsilon_{pq}\epsilon_{2pq}}$. This implies that $(2+\sqrt2)\sqrt{\epsilon_{pq}\epsilon_{2pq}}$ is a square in $\KK_3^+$. Hence  Lemma \ref{3} yields that $\left\{\epsilon_{2}, \sqrt{\epsilon_{pq}}, \sqrt{\xi\sqrt{\epsilon_{pq}\epsilon_{2pq}}}\right\}$ is a $\mathrm{F}.\mathrm{S}.\mathrm{U}$ of $\KK_3$, where $\xi$ is an $8^{th}$ root of unity.
 \end{enumerate}
  \item If  $x\pm1$ is a square in $\NN$ and  $a+1$, $a-1$ are not, then $\left\{\epsilon_{2}, \epsilon_{pq}, \sqrt{\epsilon_{2pq}}\right\}$ is a $\mathrm{F}.\mathrm{S}.\mathrm{U}$ of $\KK_3^+$  and, by  Remark \ref{37}, of  $\KK_3$.
  \item If  $a\pm1$ is a square in $\NN$ and  $x+1$, $x-1$ are not, then  $\left\{\epsilon_{2}, \epsilon_{2pq}, \sqrt{\epsilon_{pq}}\right\}$ is a $\mathrm{F}.\mathrm{S}.\mathrm{U}$ of $\KK_3^+$ and, by  Remark \ref{37}, of $\KK_3$.
  \item If $x+1$, $x-1$, $a+1$ and $a-1$ are not squares in $\NN$, then  $\left\{\epsilon_{2}, \epsilon_{pq}, \sqrt{\epsilon_{pq}\epsilon_{2pq}}\right\}$ is a  $\mathrm{F}.\mathrm{S}.\mathrm{U}$ of  $\KK_3^+$ and, by  Remark  \ref{37}, of $\KK_3$.
\end{enumerate}
\end{proof}
\vspace*{-0.5cm}
\section{\textbf{The  ambiguous classes of $\kk/\QQ(i)$}}
 Let $F=\QQ(i)$ and $\kk=\QQ(\sqrt{2pq}, i)$. We denote by  $\mathrm{A}m(\kk/F)$  the group of the ambiguous classes of $\kk/F$ and by $\mathrm{A}m_s(\kk/F)$ the subgroup of $\mathrm{A}m(\kk/F)$ generated by  the strongly ambiguous classes.
As  $p\equiv1\pmod4$, so there exist $e$ and $f$ in $\NN$ such that  $p=e^2+4f^2=\pi_1\pi_2$. Put  $\pi_1=e+2if$ and $\pi_2=e-2if$. Let  $\mathcal{H}_j$ (resp. $\o$) be the prime ideal of  $\kk$ above $\pi_j$ (resp. $1+i$), $j\in\{1, 2\}$. It is easy to see  that $\mathcal{H}_j^2=(\pi_j)$ (resp. $\o^2=(1+i)$). Therefore $[\mathcal{H}_j]\in \mathrm{A}m_s(\kk/F)$, for all  $j\in\{0, 1, 2\}$.  Keep the notation $\epsilon_{2pq}=x+y\sqrt{2pq}$. In this section,  we will determine  generators of $\mathrm{A}m_s(\kk/F)$ and $\mathrm{A}m(\kk/F)$.  Let us first prove the following result.
\begin{lem}\label{7}
Consider the prime ideals $\mathcal{H}_j$ of $\kk$, $0\leq j\leq2$.
\begin{enumerate}[\rm\indent(1)]
  \item If $x\pm1$ is  a square in $\NN$, then  $\left|\langle[\o], [\h], [\hh]\rangle\right|=8$.
  \item Else, $[\h]=[\hh]$ and $\left|\langle[\o], [\h]\rangle\right|=4$.
\end{enumerate}
\end{lem}
\begin{proof}
Since $\o^2=(1+i)$, $\mathcal{H}_l^2=(\pi_l)$ and $(\o\mathcal{H}_l)^2=((1+i)\pi_l)=\left((e\mp2f)+i(e\pm2f)\right)$,  where $1\leq l\leq2$, and since also $\sqrt{2}\not\in\QQ(\sqrt{2pq})$, $\sqrt{e^2+(2f)^2}=\sqrt{p}\not\in\QQ(\sqrt{2pq})$ and $\sqrt{(e\mp2f)^2+(e\pm2f)^2}=\sqrt{2p}\not\in\QQ(\sqrt{2pq})$, so   according to \cite[Proposition 1]{AZT12-2},  $\mathcal{H}_0$, $\mathcal{H}_l$ and $\mathcal{H}_0\mathcal{H}_l$ are not  principal  in $\kk$.\par
(1) If $x\pm1$ is  a square in $\NN$, then $p(x+1)$,  $p(x-1)$,  $2p(x+1)$ and  $2p(x-1)$ are not squares in $\NN$. Moreover $(\h\hh)^2=(p)$, hence according to \cite[Proposition 2]{AZT12-2},  $\h\hh$ is not principal in $\kk$, and the result derived.\par
(2) If $x+1$ and $x-1$ are not squares  in $\NN$, then $p(x\pm1)$ or  $2p(x\pm1)$  is a square in $\NN$; as $(\h\hh)^2=(p)$, hence according to \cite[Proposition 2]{AZT12-2},  $\h\hh$ is  principal in $\kk$. This completes the proof.
\end{proof}

 Determine now  generators of $ \mathrm{A}m_s(\kk/F)$ and $\mathrm{A}m(\kk/F)$. According to the ambiguous class number formula (see \cite{Ch-33}),
the genus number, $[(\kk/F)^*:\kk]$, is given by:
 \begin{equation}\label{51}
|\mathrm{A}m(\kk/F)|=[(\kk/F)^*:\kk]=\frac{h(F)2^{t-1}}{[E_F: E_F\cap N_{\kk/F}(\kk^\times)]},
\end{equation}
where $h(F)$ is the class number of $F$ and $t$ is the number of finite and infinite primes of $F$ ramified in $\kk/F$. Moreover as the class number
of  $F$ is equal to $1$, so the formula \eqref{51} yields that
 \begin{equation}\label{56}|\mathrm{A}m(\kk/F)|=[(\kk/F)^*:\kk]=2^r,\end{equation}
 where $r=\text{rank}\mathbf{C}l_2(\mathds{k})=t-e-1$ and $2^e=[E_F: E_F\cap N_{\kk/F}(\kk^\times)]$ (see for example \cite{McPaRa-95}).
  The relation between  $|\mathrm{A}m(\kk/F)|$ and $|\mathrm{A}m_s(\kk/F)|$ is given by the following formula (see for example \cite{Lem-13}):
\begin{equation}\label{50}
\frac{|\mathrm{A}m(\kk/F)|}{|\mathrm{A}m_s(\kk/F)|}=[E_F\cap N_{\kk/F}(\kk^\times):N_{\kk/F}(E_\kk)].
\end{equation}
 To continue, we need  the following lemma.
\begin{lem}\label{57}
Let $p\equiv -q\equiv1\pmod4$ be different primes,  $F=\QQ(i)$ and $\kk=\QQ(\sqrt{2pq}, i)$.
\begin{enumerate}[\rm\indent(1)]
  \item If $p\equiv1 \pmod8$, then $i$ is  a norm in $\kk/F$.
  \item If $p\equiv5 \pmod8$, then $i$ is not a norm in $\kk/F$.
\end{enumerate}
\end{lem}
\begin{proof}
Let  $\mathfrak{p}$ be a prime ideal  of  $F=\QQ(i)$ such that $\mathfrak{p}\neq\mathfrak{2}_F$, where $\mathfrak{2}_F$ is the prime ideal of $F$ above $2$,  then the Hilbert symbol yields that $\displaystyle\left(\frac{2pq,i}{\mathfrak{p}}\right)=\displaystyle\left(\frac{pq, i}{\mathfrak{p}}\right)$, since $2i=(1+i)^2$. Hence, by Hilbert symbol properties and  according to \cite[p. 205]{Gr-03},  we get:
\begin{itemize}
  \item If $\mathfrak{p}$ is not above $p$ and $q$, then $v_\mathfrak{p}(pq)=0$, thus
   $\displaystyle\left(\frac{pq, i}{\mathfrak{p}}\right)=1$.
  \item If $\mathfrak{p}$ lies above $p$,  then $v_\mathfrak{p}(pq)=1$, so\\
   $ \displaystyle\left(\frac{pq, i}{\mathfrak{p}}\right)=
      \displaystyle\left(\frac{ i}{\mathfrak{p}}\right)=\left(\frac{2}{p}\right),$
       indeed $\left(\frac{2}{p}\right)\left(\frac{i}{\mathfrak{p}}\right)= \left(\frac{2}{\mathfrak{p}}\right)\left(\frac{i}{\mathfrak{p}}\right)=
 \left(\frac{2i}{\mathfrak{p}}\right)=1.$
  \item If $\mathfrak{p}$ lies above $q$,  then $v_\mathfrak{p}(pq)=1$, so\\
   $ \displaystyle\left(\frac{pq, i}{\mathfrak{p}}\right)=
      \displaystyle\left(\frac{ i}{\mathfrak{p}}\right)=\left(\frac{N_{F/\QQ}(i)}{q}\right)=\left(\frac{1}{q}\right)=1,$ since $q$ remained inert in $F/\QQ$.
\end{itemize}
 So for every prime ideal  $\mathfrak{p}\in F$ and by  the product formula for the Hilbert symbol, we deduce that  $\displaystyle\left(\frac{pq, i}{\mathfrak{p}}\right)=1$, hence:\\
$(1)$ If $p\equiv1 \pmod8$, then $i$ is a norm in $\kk/F$.\\
$(2)$ If $p\equiv5\pmod8$, then $i$ is not a norm in $\kk/F$.
\end{proof}
\begin{propo}\label{248}
Let $(\kk/F)^*$ denote the relative genus field of $\kk/F$.
\begin{enumerate}[\rm\indent(1)]
\item $\k\subseteq (\kk/F)^*$ and $[(\kk/F)^*:\k]\leq2$.
\item Assume  $p\equiv1\pmod8$.
\begin{enumerate}[\rm(i)]
  \item If $x\pm1$ is a square in  $\NN$, then $\mathrm{Am}(\kk/\QQ(i))=\mathrm{Am}_s(\kk/\QQ(i))=\langle[\o], [\h], [\hh]\rangle.$
  \item Else, there  exist an unambiguous ideal $\mathcal{I}$ in $\kk/\QQ(i)$ of order  $2$ such that  $\mathrm{Am}_s(\kk/\QQ(i))=\langle[\o], [\h]\rangle$ and  $\mathrm{Am}(\kk/\QQ(i))=\langle[\o], [\h], [\mathcal{I}]\rangle.$
\end{enumerate}
\item Assume  $p\equiv5\pmod8$, then neither $x + 1$ nor $x - 1$ is a square in $\NN$ and $\mathrm{Am}(\kk/\QQ(i))=\mathrm{Am}_s(\kk/\QQ(i))=\langle[\o], [\h]\rangle$.
\end{enumerate}
\end{propo}
\begin{proof}
(1) As $\kk=\QQ(\sqrt{2pq}, i)$, so $[\k:\kk]=4$. Moreover, according to \cite[Proposition 2, p. 90]{McPaRa-95}, $r=\text{rank}\mathbf{C}l_2(\mathds{k})=3$ if $p\equiv1 \pmod8$ and $r=\text{rank}\mathbf{C}l_2(\mathds{k})=2$ if $p\equiv5 \pmod8$, so $[(\kk/F)^*:\kk]=4 \text{ or } 8$. Hence $[(\kk/F)^*:\k]\leq2$, and the   result derived.\par
(2) Assume that  $p\equiv1\pmod8$, hence $i$ is a norm in $\kk/\QQ(i)$, thus Formula \eqref{50} yields that
    \begin{align*}\dfrac{|\mathrm{Am}(\kk/\QQ(i))|}{|\mathrm{Am}_s(\kk/\QQ(i))|}&=[E_{\QQ(i)}\cap N_{\kk/\QQ(i)}(\kk^{\times}):N_{\kk/\QQ(i)}(E_\kk)]\\
                                                                               &=\left\{\begin{array}{ll}
                                                                                 1 \text{ if  $x\pm1$ is a square in } \NN,\\
                                                                                 2 \text{ if not, }
                                                                                 \end{array}\right.\end{align*}
since in the case where $x\pm1$ is a square in  $\NN$, we have $E_\kk=\langle i, \sqrt{i\epsilon_{2pq}}\rangle$, hence $[E_{\QQ(i)}\cap
N_{\kk/\QQ(i)}(\kk^{\times}):N_{\kk/\QQ(i)}(E_\kk)]=[< i >: < i >]=1$, and if not we have  $E_\kk=\langle i, \epsilon_{2pq}\rangle$, hence $[E_{\QQ(i)}\cap
N_{\kk/\QQ(i)}(\kk^{\times}):N_{\kk/\QQ(i)}(E_\kk)]=[< i >: < -1>]=2$.\\
\indent On the other hand, as  $p\equiv1\pmod8$, so   according to \cite[Proposition 2, p. 90]{McPaRa-95},  $r=3$. Therefore     $|\mathrm{Am}(\kk/\QQ(i))|=2^3.$\\
\indent (i) If $x\pm1$ is a square in $\NN$, then $\mathrm{Am}_s(\kk/\QQ(i))=\mathrm{Am}(\kk/\QQ(i))$, hence by Lemma \ref{7} we get  $\mathrm{Am}(\kk/\QQ(i))=\mathrm{Am}_s(\kk/\QQ(i))=\langle[\o], [\h], [\hh]\rangle.$\\
\indent (ii)  If $x+1$ and $x-1$ are not squares in $\NN$, then $$|\mathrm{Am}(\kk/\QQ(i))|=2|\mathrm{Am}_s(\kk/\QQ(i))|=8,$$  hence Lemma \ref{7} yields that $\mathrm{Am}_s(\kk/\QQ(i))=\langle[\o], [\h]\rangle.$\\
\indent Consequently,  there exist an unambiguous  ideal $\mathcal{I}$ in $\kk/F$ of order $2$ such that
    $$\mathrm{Am}(\kk/\QQ(i))=\langle[\o], [\h], [\mathcal{I}]\rangle.$$
    By Chebotarev theorem,  $\mathcal{I}$ can  always be chosen  as a prime ideal of $\kk$ above a prime $l$ in $\QQ$, which splits completely in $\kk$.
    So we can determine  $\mathcal{I}$ by using the following lemma.
    \begin{lem}[\cite{Si-95}]
     Let $p_1$, $p_2$,...,$p_n$ be distinct primes and for each $j$,  let $e_j=\pm1$.
  Then there exist infinitely many primes $l$ such that $\left(\frac{p_j}{l}\right)=e_j$, for all $j$.
    \end{lem}
    Let $l\equiv1\pmod4$ be a prime satisfying   $\left(\frac{2pq}{l}\right)=-\left(\frac{q}{l}\right)=1$, then $l$ splits completely in  $\kk$. Let
    $\mathcal{I}$ be a prime ideal of  $\kk$  above $l$;  hence  $\mathcal{I}$ remained inert in  $\KK_2$ and $\left(\frac{2p}{l}\right)=-1$. We need to prove that $\mathcal{I}$,  $\o\mathcal{I}$,  $\h\mathcal{I}$ and $\o\h\mathcal{I}$ are not principal in  $\kk$.\\
     \indent$\bullet$ As $\mathcal{I}$ remained inert in $\KK_2$, so $\varphi_{\KK_2/\kk}(\mathcal{I})\neq1$, where $\varphi_{\KK_2/\kk}$ denotes the Artin map of $\KK_2$ over $\kk$, similarly, we have  $\varphi_{\KK_2/\kk}(\h\mathcal{I})\neq1$ (note that $\left(\frac{p}{q}\right)=1$, since $p(x\pm1)$ or $2p(x\pm1)$ is a square in  $\NN$). Therefore $\mathcal{I}$ and $\h\mathcal{I}$ are not principal in  $\kk$.\\
      \indent$\bullet$ Let us prove that $\o\mathcal{I}$ is not  principal in $\kk$.  For this, we consider the following cases:
      \begin{enumerate}[\indent(a)]
        \item Assume $\left(\frac{2}{l}\right)=1$, then $\left(\frac{p}{l}\right)=-1$;  thus
       if $\left(\frac{2}{q}\right)=-1$, then $\varphi_{\KK_3/\kk}(\o\mathcal{I})\neq1$, and if $\left(\frac{2}{q}\right)=1$, then $\varphi_{\KK_1/\kk}(\o\mathcal{I})\neq1$. Hence $\o\mathcal{I}$ is not principal in $\kk$.
        \item Assume now $\left(\frac{2}{l}\right)=-1$, hence $\left(\frac{p}{l}\right)=1$. Thus
     if $\left(\frac{2}{q}\right)=1$, then $\varphi_{\KK_2/\kk}(\o\mathcal{I})\neq1$. If
      $\left(\frac{2}{q}\right)=-1$, so we need   the following two  quadratic extensions of $\kk$:   $\KK_4=\kk(\sqrt{\pi_1})$ and $\KK_5=\kk(\sqrt{2\pi_1})=\kk(\sqrt{\pi_2q})$, where $p=e^2+16f^2=\pi_1\pi_2=(e+4if)(e-4if)$, since $p\equiv1\pmod8$. Note that  $\KK_4/\kk$ and $\KK_5/\kk$ are unramified (see \cite{AZT12-1}).  As $\left(\frac{2}{p}\right)=1$, then $\left(\frac{1+i}{\pi_1}\right)=\left(\frac{1+i}{\pi_2}\right)$, hence the quadratic residue symbol implies that $$\left(\frac{\pi_1}{\o\mathcal{I}}\right)=\left(\frac{1+i}{\pi_1}\right)=-\left(\frac{\pi_2q}{\o\mathcal{I}}\right).$$
      Therefore, if $\left(\frac{1+i}{\pi_1}\right)=-1$ , then $\varphi_{\KK_4/\kk}(\o\mathcal{I})\neq1$, else we have $\varphi_{\KK_5/\kk}(\o\mathcal{I})\neq1$. Thus  $\o\mathcal{I}$ is not  principal in $\kk$.
      \end{enumerate}

With the same argument, we show that  $\o\h\mathcal{I}$ is not principal in $\kk$.

(3) Assume that  $p\equiv5\pmod8$, hence $i$ is not a norm in $\kk/\QQ(i)$ and $x+1$, $x-1$ are not squares in $\NN$, for if $x\pm1$ is a square in $\NN$, then the Legendre symbol implies that $$1=\left(\frac{x\pm1}{p}\right)=\left(\frac{x\mp1\pm2}{p}\right)=\left(\frac{2}{p}\right),$$ which is absurd.  Thus $|\mathrm{Am}(\kk/\QQ(i))|=2^2$ and
    $$\dfrac{|\mathrm{Am}(\kk/\QQ(i))|}{|\mathrm{Am}_s(\kk/\QQ(i))|}=[E_{\QQ(i)}\cap N_{\kk/\QQ(i)}(\kk^{\times}):N_{\kk/\QQ(i)}(E_\kk)]=1.$$                                                                         Hence by Lemma \ref{7} we get   $\mathrm{Am}(\kk/\QQ(i))=\mathrm{Am}_s(\kk/\QQ(i))=\langle[\o], [\h]\rangle$.
      This completes the proof.
\end{proof}
\section{\bf{Capitulation}}
In this section,  we will determine the classes of  $\mathbf{C}l_2(\kk)$, the  $2$-class group of $\kk$, that capitulate in  $\KK_j$, for all $j\in\{1, 2, 3\}$. For this we need the following theorem.
\begin{The}[\cite{HS82}]\label{1}
 Let $K/k$ be a cyclic extension of prime degree, then  the number of classes that capitulate in $K/k$ is:
 $[K:k][E_k:N_{K/k}(E_K)],$
 where $E_k$ and $E_K$ are the  unit groups of $k$ and $K$ respectively.
\end{The}
\begin{The}\label{226}
Let $\KK_j$, $1\leq j\leq3$, be the three unramified quadratic extensions of $\kk$ defined above.
\begin{enumerate}[\rm\indent(1)]
\item For $j\in\{1, 2\}$ we have
\begin{enumerate}[\rm(i)]
  \item If  $x\pm1$ is a square in $\NN$, then $|\kappa_{\KK_j}|=4$.
  \item Else $|\kappa_{\KK_j}|=2$.
\end{enumerate}
\item Put $\epsilon_{pq}=a+b\sqrt{pq}$ and let $Q_{\KK_3}$ denote the unit index of  $\KK_3$.
\begin{enumerate}[\rm(i)]
\item If  both of $x\pm1$ and  $a\pm1$  are squares in $\NN$, then
\begin{enumerate}[\rm(a)]
  \item If  $Q_{\KK_3}=1$, then $|\kappa_{\KK_3}|=4$.
  \item If $Q_{\KK_3}=2$, then $|\kappa_{\KK_3}|=2$.
\end{enumerate}
\item If one of the four numbers   $x+1$,  $x-1$,  $a+1$ and $a-1$ is a square in  $\NN$, and the others are not, then $|\kappa_{\KK_3}|=4$.
\item If $x+1$, $x-1$,  $a+1$ and $a-1$ are not squares in $\NN$, then $|\kappa_{\KK_3}|=2$.
\end{enumerate}
\end{enumerate}
\end{The}
\begin{proof}
$(1)$ According to   Proposition $\ref{27}$,  $E_{\KK_1}=\langle i, \epsilon_{p}, \sqrt{i\epsilon_{2q}}, \sqrt{\epsilon_{2q}\epsilon_{2pq}}\rangle$ or\\ $\langle i, \epsilon_{p}, \sqrt{i\epsilon_{2q}}, \sqrt{\epsilon_{2pq}}\rangle$, so $N_{\KK_1/\kk}(E_{\KK_1})=\langle i, \epsilon_{2pq}\rangle$. On the other hand, Proposition \ref{28} yields that $E_{\KK_2}=\langle i, \sqrt{i\epsilon_{q}}, \sqrt{i\epsilon_{2p}}, \sqrt{i\epsilon_{2pq}}\rangle$ or
  $\langle i, \sqrt{i\epsilon_{q}}, \sqrt{i\epsilon_{2p}}, \sqrt{\epsilon_{2pq}}\rangle$ or\\ $\langle i, \sqrt{i\epsilon_{q}},
   \epsilon_{2p}, \sqrt{\epsilon_{q}\epsilon_{2pq}}\rangle$  or $\langle i, \sqrt{i\epsilon_{q}}, \epsilon_{2p},
   \sqrt{\epsilon_{2pq}}\rangle$, hence $N_{\KK_2/\kk}(E_{\KK_2})=\langle i, \epsilon_{2pq}\rangle$.\\
 \indent (i) If  $x\pm1$ is a square in $\NN$, then   Lemma \ref{6} yields that  $E_{\kk}=\langle i, \sqrt{i\epsilon_{2pq}}\rangle$. Therefore  $[E_{\kk}:N_{\KK_j/\kk}(E_{\KK_j})]=2$, and  Theorem \ref{1} implies that    $|\kappa_{\KK_j}|=4$.\\
 \indent (ii) Else $E_{\kk}=\langle i, \epsilon_{2pq}\rangle$, which gives that $[E_{\kk}:N_{\KK_j/\kk}(E_{\KK_j})]=1$, and  Theorem \ref{1} implies that  $|\kappa_{\KK_1}|=2$.\\
 \indent (2)  (i) Assume  that  $x\pm1$ and  $a\pm1$  are squares in $\NN$, so by  Lemma \ref{6} we get $E_{\kk}=\langle i, \sqrt{i\epsilon_{2pq}}\rangle$.
\begin{enumerate}[\rm(a)]
 \item If $Q_{\KK_3}=1$, then Proposition $\ref{29}$ implies that  $E_{\KK_3}=\langle \sqrt i, \epsilon_{2}, \sqrt{\epsilon_{pq}}, \sqrt{\epsilon_{2pq}}\rangle$, hence $N_{\KK_3/\kk}(E_{\KK_3})=\langle i, \epsilon_{2pq}\rangle$,  from which we deduce that $[E_{\kk}:N_{\KK_3/\kk}(E_{\KK_3})]=2$, and  Theorem \ref{1} implies that   $|\kappa_{\KK_3}|=4$.
 \item If $Q_{\KK_3}=2$, then  Proposition $\ref{29}$ implies that  $E_{\KK_3}=\langle \sqrt i, \epsilon_{2}, \sqrt{\epsilon_{pq}}, \sqrt{\sqrt{i\epsilon_{pq}\epsilon_{2pq}}}\rangle$, thus $N_{\KK_3/\kk}(E_{\KK_3})=\langle i, \sqrt{i\epsilon_{2pq}}\rangle$, from which we deduce that $[E_{\kk}:N_{\KK_3/\kk}(E_{\KK_3})]=1$,  and  Theorem \ref{1} implies that $|\kappa_{\KK_3}|=2$.
 \end{enumerate}

 \indent (ii) If $x\pm1$ is a square in $\NN$ and $a+1$,  $a-1$ are not, then by Lemma \ref{6} we get $E_{\kk}=\langle i, \sqrt{i\epsilon_{2pq}}\rangle$. Moreover,   Proposition $\ref{29}$ implies that  $E_{\KK_3}=\langle \sqrt i, \epsilon_{2}, \epsilon_{pq}, \sqrt{\epsilon_{2pq}}\rangle$, hence $N_{\KK_3/\kk}(E_{\KK_3})=\langle i, \epsilon_{2pq}\rangle$. Therefore  $[E_{\kk}:N_{\KK_3/\kk}(E_{\KK_3})]=2$,   and  Theorem \ref{1} implies that   $|\kappa_{\KK_3}|=4$.\\
 If $a\pm1$ is a square in $\NN$ and $x+1$, $x-1$ are not, then by Lemma \ref{6} we get $E_{\kk}=\langle i, \epsilon_{2pq}\rangle$. Moreover  Proposition $\ref{29}$ implies that  $E_{\KK_3}=\langle \sqrt i, \epsilon_{2}, \sqrt{\epsilon_{pq}}, \epsilon_{2pq}\rangle$, hence  $N_{\KK_3/\kk}(E_{\KK_3})=\langle i, \epsilon_{2pq}^2\rangle$. Therefore  $[E_{\kk}:N_{\KK_3/\kk}(E_{\KK_3})]=2$,  and  Theorem \ref{1} implies that  $|\kappa_{\KK_3}|=4$.\\
 \indent (iii) Finally, assume that $x+1$, $x-1$, $a+1$ and $a-1$ are not squares in $\NN$, then by  Lemma \ref{6} we get $E_{\kk}=\langle i, \epsilon_{2pq}\rangle$. Moreover,  Proposition $\ref{29}$ implies that  $E_{\KK_3}=\langle \sqrt i, \epsilon_{2}, \epsilon_{pq}, \sqrt{\epsilon_{pq}\epsilon_{2pq}}\rangle$, hence $N_{\KK_3/\kk}(E_{\KK_3})=\langle i, \epsilon_{2pq}\rangle$. Therefore  $[E_{\kk}:N_{\KK_3/\kk}(E_{\KK_3})]=1$,  and  Theorem \ref{1} implies that  $|\kappa_{\KK_3}|=2$.
\end{proof}
\subsection{Capitulation in $\KK_1$}
\begin{The}\label{227}
Keep  the notations and hypotheses previously  mentioned.
\begin{enumerate}[\rm\indent(1)]
 \item If  $x\pm1$ is a square in $\NN$, then $\kappa_{\KK_1}=\langle [\mathcal{H}_1], [\mathcal{H}_2]\rangle$.
  \item Else $\kappa_{\KK_1}=\langle[\mathcal{H}_1]\rangle$.
\end{enumerate}
\end{The}
\begin{proof}
 Let us first  prove  that  $\h$ and $\hh$ capitulate in  $\KK_1$.
As $N(\epsilon_{p})=-1$, then $s^2+4=t^2p$, where  $\epsilon_{p}=\frac{1}{2}(s+t\sqrt{p})$, hence
$(s-2i)(s+2i)=t^2p$.  According to the decomposition uniqueness in $\ZZ[i]$, there exist  $t_1$ and $t_2$ in $\ZZ[i]$ such that:
 $(1) \left\{\begin{array}{ll}
s\pm2i &=t_1^2\pi_1\\
s\mp2i &= t_2^2\pi_2,
\end{array}\right. \text{ or }
 (2) \left\{\begin{array}{ll}
s\pm2i &=it_1^2\pi_1\\
s\mp2i &=-it_2^2\pi_2,
\end{array}\right.\text{ where } t=t_1t_2.$

\indent $\bullet$ The  system (1) implies that  $2s=t_1^2\pi_1+t_2^2\pi_2$.   Put
$\alpha = \frac{1}{2}(t_1\pi_1+t_2\sqrt p)$ and   $\beta =\frac{1}{2}(
t_2\pi_2+t_1\sqrt p)$. Then  $\alpha$ and $\beta$ are in  $\KK_1=\kk(\sqrt
p)$ and we have:
 \begin{align*}
  \alpha^2 & = \frac{1}{4}(t_1^2\pi_1^2+t_2^2p+2t_1t_2\pi_1\sqrt p)\\
           &= \frac{1}{4}\pi_1(t_1^2\pi_1+t_2^2\pi_2+2t\sqrt p),\ \text{ since }\  p=\pi_1\pi_2\  and\ t=t_1t_2.\\
           &= \frac{1}{4}\pi_1(2s+2t\sqrt p),\ \text{ since }\  2s=t_1^2\pi_1+t_2^2\pi_2.\\
           &= \pi_1\epsilon_{p},\  \text{ since }\  \epsilon_{p}=\frac{1}{2}(s+t\sqrt p).
  \end{align*}
  The same argument yields that
$\beta^2 = \pi_2\epsilon_{p}$.\\
\indent Consequently,
   $(\alpha^2)=(\pi_1)=\mathcal{H}_1^2$ (resp.  $(\beta^2)=(\pi_2)=\mathcal{H}_2^2$),
   hence $(\alpha)= \mathcal{H}_1$ and $(\beta)= \mathcal{H}_2$. \\
\indent  $\bullet$ Similarly,   system (2) yields that  $2s=it_1^2\pi_2-it_2^2\pi_1$, hence $\sqrt{2\pi_1\epsilon_{p}}=\frac{1}{2}(t_1(1+i)\pi_1+t_2(1-i)\sqrt{p})$ and $\sqrt{2\pi_2\epsilon_{p}}=\frac{1}{2}(t_1(1+i)\sqrt{p}+t_2(1-i)\pi_2)$  are in $\KK_1$. Therefore  there exist   $\alpha$ and $\beta$ in $\KK_1$ such that $2\pi_1\epsilon_{p}=\alpha^2$ and $2\pi_2\epsilon_{p}=\beta^2$, thus $(\frac{\alpha}{1+i})= \mathcal{H}_1$ and $(\frac{\beta}{1+i})= \mathcal{H}_2$. This yields that $\h$ and $\hh$ capitulate in $\KK_1$.\\
 \indent On the other hand, by Lemma \ref{7},  $\mathcal{H}_j$, $1\leq j\leq2$,  are not  principal in $\kk$.\\
\indent  (1) If $x\pm1$ is a square in $\NN$, then  Lemma \ref{7} yields that  $[\h\hh]\neq1$. Hence the result.\\
 \indent (2) If $x+1$ and $x-1$ are not squares in $\NN$, then   Lemma \ref{7} yields that  $[\h]=[\hh]$.  This completes the proof.
\end{proof}
\subsection{Capitulation in $\KK_2$}
We need the  following two lemmas.
\begin{lem}\label{225}
If  $N(\epsilon_{2p})=1$, then
\begin{enumerate}[\rm\indent(1)]
  \item\label{item1} $p\equiv1\pmod8$.
  \item $2p(x-1)$ is not a square in $\NN$.
\end{enumerate}
\end{lem}
\begin{proof}
$(1)$ Put $\epsilon_{2p}=\alpha+\beta\sqrt{2p}$, then, if $N(\epsilon_{2p})=1$,  Lemma \ref{4} yields that
 $\left\{\begin{array}{ll}
 \alpha\pm1=\beta_1^2,\\
 \alpha\mp1=2p\beta_2^2,
 \end{array}\right.$
 hence $1=\left(\frac{\alpha\pm1}{p}\right)=\left(\frac{\alpha\mp1\pm2}{p}\right)=\left(\frac{2}{p}\right)$, so the result.\\
 $(2)$ If $2p(x-1)$ is  a square in $\NN$, then
  $\left\{\begin{array}{ll}
 x-1=2py_1^2,\\
 x+1=qy_2^2;
 \end{array}\right.$\\
 thus
 $$\left\{\begin{array}{ll}
 \left(\frac{2p}{q}\right)=\left(\frac{x-1}{q}\right)=-\left(\frac{2}{q}\right),\\
 \left(\frac{q}{p}\right)=\left(\frac{x+1}{p}\right)=\left(\frac{2}{p}\right);
 \end{array}\right.$$
 this implies that $\left(\frac{2}{p}\right)=-1$, which contradicts the first assertion $(1)$.
\end{proof}
\begin{lem}\label{247}
Put $\epsilon_{pq}=a+b\sqrt{pq}$.
If  $a\pm1$ is a square in  $\NN$,  then
 $p\equiv1\pmod8$.
\end{lem}
\begin{proof}
The same argument as in Lemma \ref{225}(\ref{item1}) leads to the result.
\end{proof}
\begin{The}\label{229}
Keep  the notations and hypotheses previously  mentioned.
\begin{enumerate}[\rm\indent(1)]
 \item If $N(\epsilon_{2p})=1$ and $x\pm1$ is a square in $\NN$, then $\kappa_{\KK_2}=\langle [\o], [\h\hh]\rangle$ or
  $\langle [\h], [\hh]\rangle$.
\item\label{item2}  If $N(\epsilon_{2p})=1$ and $x+1$,  $x-1$ are not squares in $\NN$, then there exist an unambiguous  ideal $\mathcal{I}$ in $\kk/F$ of order $2$, such that $\kappa_{\KK_2}=\langle[\mathcal{I}]\rangle$ or $\langle[\o\mathcal{I}]\rangle$ or $\langle[\h\mathcal{I}]\rangle$   or $\langle[\o\h\mathcal{I}]\rangle$.
\item  If $N(\epsilon_{2p})=-1$, then
\begin{enumerate}[\rm(i)]
 \item If  $x\pm1$ is a square in $\NN$, then $\kappa_{\KK_2}=\langle [\o\mathcal{H}_1], [\o\mathcal{H}_2]\rangle$
  \item Else, $\kappa_{\KK_2}=\langle[\o\mathcal{H}_1]\rangle$.
\end{enumerate}
\end{enumerate}
\end{The}
\begin{proof}
Since   $(\pi_j)=\mathcal{H}_j^2$,  $j\in\{1, 2\}$, and   $\o^2=(1+i)$, so $(2p)=((1+i)\h\hh)^2$. Moreover,  $2p$ is a square in $\KK_2$, so there exist  $\alpha\in\KK_2$ such that $(2p)=(\alpha^2)$,  hence $((1+i)\h\hh)^2=(\alpha^2)$, therefore $\h\hh=\left(\frac{\alpha}{1+i}\right)$ and  $\h\hh$ capitulates in $\KK_2$.\\
 \indent (1) If $N(\epsilon_{2p})=1$, then by Lemma  \ref{225} we get $p\equiv1\pmod8$. Moreover, according to Lemma \ref{7}, if  $x\pm1$ is a square in $\NN$, then $\h$, $\hh$ and $\h\hh$ are not principal in   $\kk$ and, according to the  Theorem \ref{226}, there are four classes that capitulate in $\KK_2$. The following examples affirm the two cases of capitulation:

{ \scriptsize
  \begin{longtable}{| p{0.693in} | p{0.48in} | p{0.4in} |p{0.4in} |p{0.4in} | p{0.46in} | p{0.4in} |p{0.5in} |}
\hline
 $d$  $ = 2pq$ & $x+1$ & $\o\mathcal{O}_{\KK_2}$ & $\h\mathcal{O}_{\KK_2}$ & $\hh\mathcal{O}_{\KK_2}$  & $\h\hh\mathcal{O}_{\KK_2}$ & $\mathbf{C}l(\kk)$ & $\mathbf{C}l(\KK_2)$\\
\hline
\endfirsthead
\hline
 $d$  $ = 2pq$ & $x+1$ & $\o\mathcal{O}_{\KK_2}$ & $\h\mathcal{O}_{\KK_2}$ & $\hh\mathcal{O}_{\KK_2}$  & $\h\hh\mathcal{O}_{\KK_2}$ & $\mathbf{C}l(\kk)$ & $\mathbf{C}l(\KK_2)$\\
\hline
\endhead
$238 = 2.17.7$ & $108^2$ & $[0, 0, 0]~$ & $[4, 0, 0]~$ & $[4, 0, 0]~$  & $[0, 0, 0]$ & $(4, 2, 2)$ & $(8, 2, 2)$\\ \hline
$782 = 2.17.23$ & $28^2$ & $[0, 0, 0]~$ & $[12, 0, 0]~$ & $[12, 0, 0]~$  & $[0, 0, 0]~$ & $(12, 2, 2)$ & $(24, 6, 2)$\\ \hline
$1022 = 2.73.7$ & $32^2$ & $[16, 0, 0]~$ & $[0, 0, 0]~$ & $[0, 0, 0]~ $  & $[0, 0, 0]~$ & $(16, 2, 2)$ & $(32, 8, 2)$\\ \hline
$1246 = 2.89.7$ & $21068856^2$ & $[8, 0, 0]~$ & $[0, 0, 0]~$ & $[0, 0, 0]~$  & $[0, 0, 0]~$ & $(8, 2, 2)$ & $(16, 4, 2)$\\ \hline
$1358 = 2.97.7$ & $1732^2$ & $[0, 0, 0]~$ & $[60, 0, 0]~$ & $[60, 0, 0]~$  & $[0, 0, 0]~$ & $(12, 2, 2)$ & $(120, 2, 2)$\\ \hline
\end{longtable}

\begin{longtable}{| c | c | c | c | c |c | c | c |}
\hline
 $d$  $ = 2pq$  & $x-1$ & $\o\mathcal{O}_{\KK_2}$ & $\h\mathcal{O}_{\KK_2}$ & $\hh\mathcal{O}_{\KK_2}$ & $\h\hh\mathcal{O}_{\KK_2}$ & $\mathbf{C}l(\kk)$ & $\mathbf{C}l(\KK_2)$\\
\hline
\endfirsthead
\hline
 $d$  $ = 2pq$  & $x-1$ & $\o\mathcal{O}_{\KK_2}$ & $\h\mathcal{O}_{\KK_2}$ & $\hh\mathcal{O}_{\KK_2}$ & $\h\hh\mathcal{O}_{\KK_2}$ & $\mathbf{C}l(\kk)$ & $\mathbf{C}l(\KK_2)$\\
\hline
\endhead
$374 = 2.17.11$ & $58^2$ & $[0, 2]~$ & $[0, 0]$ & $[0, 0]~$ & $[0, 0]~$ & $(14, 2, 2)$ & $(28, 4)$\\ \hline
$534 = 2.89.3$ & $1918^2$ & $[0, 0]~$ & $[40, 0]~$ & $[40, 0]~$ & $[0, 0]~$ & $(10, 2, 2)$ & $(80, 2)$\\ \hline
$1398 = 2.233.3$ & $2206^2$ & $[0, 0]~$ & $[40, 0]~$ & $[40, 0]~$ & $[0, 0]~$ & $(10, 2, 2)$ & $(80, 2)$\\ \hline
$2118 = 2.353.3$ & $46^2$ & $[60, 12]~$ & $[0, 0]~$ & $[0, 0]~$ & $[0, 0]~$ & $(30, 2, 2)$ & $(120, 24)$\\ \hline
$2694 = 2.449.3$ & $2095718^2$ & $[0, 6, 0]~$ & $[0, 0, 0]~$ & $[0, 0, 0]~$ & $[0, 0, 0]~$ & $(30, 2, 2)$ & $(60, 12, 3)$\\ \hline
\end{longtable}
}
\vspace*{-0.4cm}
   (2) If $N(\epsilon_{2p})=1$ and $x+1$,  $x-1$ are not squares in $\NN$, then we are in the assumptions of the Proposition \ref{248}, since $N(\epsilon_{2p})=1$ yields that $p\equiv1\pmod8$. Moreover, Lemma \ref{6} implies that  $E_{\kk}=\langle i, \epsilon_{2pq}\rangle$. \\
   \indent (i)  Assume   $2p(x+1)$ is a square in $\NN$, hence, according to the Proposition $\ref{28}$, we have  $E_{\KK_2}=\langle i, \sqrt{i\epsilon_{q}}, \sqrt{i\epsilon_{2p}}, \sqrt{i\epsilon_{2pq}}\rangle$ and, according to the Theorem \ref{226}, there are two classes that capitulate in $\KK_2$. So to prove the result, it suffices to show that $\o$,  $\h$ and $\o\h$ do not capitulate in  $\KK_2$.
 If $\o$ (resp. $\h$,  $\o\h$) capitulates in  $\KK_2$, then there exist $\alpha\in\KK_2$ such that $\o=(\alpha)$ (resp. $\h=(\alpha)$,  $\o\h=(\alpha)$), hence  $(\alpha^2)= (1+i)$ (resp. $(\alpha^2)= (\pi_1)$, $(\alpha^2)= ((1+i)\pi_1)$). Consequently, $(1+i)\epsilon=\alpha^2$ (resp. $\alpha^2= \pi_1\epsilon$, $\alpha^2= (1+i)\pi_1\epsilon$)  with some unit $\epsilon\in\KK_2$; note that $\epsilon$ can be taken as follows $\epsilon=i^a\sqrt{i\epsilon_{q}}^b\sqrt{i\epsilon_{2p}}^c\sqrt{i\epsilon_{2pq}}^d$, where $a$, $b$, $c$ and $d$ are in  $\{0, 1\}$.

  At first, let us show that the unit   $\epsilon$ is neither real nor purely imaginary. In fact, if it is real  (same proof if it is purely imaginary), then putting $\alpha=\alpha_1+i\alpha_2$, where  $\alpha_j\in\KK_2^+$, we get:
\begin{enumerate}[a.]
  \item If $(1+i)\epsilon=\alpha^2$,  then  $\alpha_1^2-\alpha_2^2+2i\alpha_1\alpha_2=\epsilon(1+i)$, hence
$\left\{
 \begin{array}{ll}
 \alpha_1^2-\alpha_2^2=\epsilon,\\
 2\alpha_1\alpha_2=\epsilon,
 \end{array}\right.$\\
  thus $\alpha_1^2-2\alpha_2\alpha_1-\alpha_2^2=0$; therefore  $\alpha_1=\alpha_2(1\pm\sqrt 2)$, and  $\sqrt 2\in\KK_2^+$ (for the case  $\alpha^2= \pi_1\epsilon$, we get $\sqrt p\in\KK_2^+$), which is absurd.
  \item If $(1+i)\pi_1\epsilon=\alpha^2$,  then  $\alpha_1^2-\alpha_2^2+2i\alpha_1\alpha_2=\epsilon(1+i)\pi_1$, hence
$\left\{
 \begin{array}{ll}
 \alpha_1^2-\alpha_2^2=\epsilon(e-4f),\\
 2\alpha_1\alpha_2=\epsilon(e+4f),
 \end{array}\right.$
 where $p=e^2+16f^2$, since $p\equiv1\pmod8$.
  Thus $$4\alpha_1^4-4\epsilon(e-4f)\alpha_1^2-\epsilon^2(e+4f)^2=0,$$
  from which we deduce that $\alpha_1^2=\frac{\epsilon}{2}\left[(e-4f)\pm\sqrt{2p}\right]$. As $\alpha_1\in\KK_2^+$, so putting $\alpha_1=a+b\sqrt{2p}$, where $a$, $b$ are in $\QQ(\sqrt q)$,  we get the following unsolvable equation (in $\QQ(\sqrt q)$):
  $$16a^4-8\epsilon(e-4f)a^2+2p\epsilon^2=0,$$ since its reduced discriminant is $\Delta'=-16\epsilon^2(e+4f)^2<0$.
\end{enumerate}
 To this end, as $(1+i)\epsilon=\alpha^2$ (same proof for the other cases), then applying the norm $N_{\KK_2/\kk}$ we get that $(1+i)^2N_{\KK_2/\kk}(\epsilon)=N_{\KK_2/\kk}(\alpha)^2$, with  $N_{\KK_2/\kk}(\epsilon)\in E_{\kk}=\langle i, \epsilon_{2pq}\rangle$. Without loss of generality, one can take $N_{\KK_2/\kk}(\epsilon)\in\{\pm1, \pm i, \pm \epsilon_{2pq}, \pm i\epsilon_{2pq}\}$.
\begin{itemize}
 \item  As $N_{\KK_2/\kk}(\epsilon)$ is a square in $E_{\kk}$, so $N_{\KK_2/\kk}(\epsilon)\not\in\{ \pm i, \pm \epsilon_{2pq}, \pm i\epsilon_{2pq}\}$.
  \item  If $N_{\KK_2/\kk}(\epsilon)=\pm1$, then there exist  $a$, $b$, $c$ and $d$ in $\{0, 1\}$ such that $\epsilon=i^a\sqrt{i\epsilon_{q}}^b\sqrt{i\epsilon_{2p}}^c\sqrt{i\epsilon_{2pq}}^d$ and $N_{\KK_2/\kk}(\epsilon)=\pm1$, hence, $(-1)^a\epsilon_{2pq}^{d}i^{b+c+d}=\pm1$; so necessarily we must have $b=c$ and $d=0$. Therefore $\epsilon=i^{a+b}\sqrt{\epsilon_{q}\epsilon_{2p}}^b$, which contradicts  the fact that $\epsilon$  is not  real or purely imaginary.
\end{itemize}
The following examples clarify this:  the first table gives  examples of the ideals $\mathcal{I}$,  $\o$ and  $\h$ which are not principal in $\kk$, and gives the structures of the class groups of $\kk$ and $\KK_2$ respectively; whereas the second table gives the cases of  capitulation  of these ideals in $\KK_2$.
{\tiny
\begin{longtable}{| p{0.84in} | p{0.5in} | p{0.4in}| p{0.4in} |p{0.4in} |p{0.4in} | p{0.5in} | p{0.5in} |}
\hline
$d = 2pq$  & $2p(x+1)$ & $\mathcal{I}$ &  $\mathcal{I}^2$ & $\o$ & $\h$ &  $\mathbf{C}l(\kk)$ & $\mathbf{C}l(\KK_2)$\\
\hline
\endfirsthead
\hline
 $d = 2pq$  & $2p(x+1)$ & $\mathcal{I}$ & $\mathcal{I}^2$ & $\o$ & $\h$  & $\mathbf{C}l(\kk)$ & $\mathbf{C}l(\KK_2)$\\
\hline
\endhead
$582 = 2.97.3$ & $194^2$ & $[0, 1, 1]$ & $[0,  0, 0]$ & $[4, 1, 1]$ & $[4, 0, 0]$ & $(8, 2, 2)$ & $(80, 4, 2)$\\ \hline
$646 = 2.17.19$ & $102^2$ & $[4, 0, 0]~$ & $[0,  0, 0]$ & $[0, 0, 1]~$ & $[4, 2, 0]~$ & $(8, 4, 2)$ & $(8, 8, 2, 2)$ \\ \hline
$2822 = 2.17.83$ & $850^2$ & $[12, 1, 0]~$ & $[0,  0, 0]$ & $[0, 0, 1]~$ & $[12, 0, 0]~$ & $(24, 2, 2)$ & $(48, 12, 2)$\\  \hline
$5654 = 2.257.11$ & $178358^2$ & $[28, 1, 0]~$ & $[0,  0, 0]$ & $[0, 0, 1]~$ & $[28, 0, 0]~$ & $(56, 2, 2)$ & $(224, 8, 4)$\\ \hline
$8854 = 2.233.19$ & $9786^2$ & $[0, 0, 1]$ & $[0,  0, 0]$ & $[60, 0, 1]~$ & $[60, 0, 0]~$ & $(120, 2, 2)$ & $(120, 8, 2, 2)$\\ \hline
$10806 = 2.1801.3$ & $258569570^2$ & $[0, 1, 1]~$ & $[0,  0, 0]$ & $[24, 1, 0]~$ & $[24, 0, 0]~$ & $(48, 2, 2)$ & $(48, 48, 6, 2)$\\ \hline
\end{longtable}
\begin{longtable}{| p{0.821in}  | p{0.475in} |p{0.5in} |p{0.43in} | p{0.43in} | p{0.5in} |p{0.45in} |p{0.5in} |}
\hline
 $d = 2pq$  & $\o\mathcal{O}_{\KK_2}$ & $\h\mathcal{O}_{\KK_2}$ & $\o\h\mathcal{O}_{\KK_2}$ & $\mathcal{I}\mathcal{O}_{\KK_2}$ & $\h\mathcal{I}\mathcal{O}_{\KK_2}$ & $\o\mathcal{I}\mathcal{O}_{\KK_2}$ & $\o\h\mathcal{I}\mathcal{O}_{\KK_2}$\\
\hline
\endfirsthead
\hline
 $d = 2pq$  & $\o\mathcal{O}_{\KK_2}$ & $\h\mathcal{O}_{\KK_2}$ & $\o\h\mathcal{O}_{\KK_2}$ & $\mathcal{I}\mathcal{O}_{\KK_2}$ & $\h\mathcal{I}\mathcal{O}_{\KK_2}$ & $\o\mathcal{I}\mathcal{O}_{\KK_2}$ & $\o\h\mathcal{I}\mathcal{O}_{\KK_2}$\\
\hline
\endhead
$582 = 2.97.3$ &    $[0, 2, 0]~$ & $[40, 2, 0]$  & $[40, 0, 0]~$ & $[40, 2, 0]~$ & $[0, 0, 0]~$ & $[40, 0, 0]$ & $[0, 2, 0]$\\ \hline
$646 = 2.17.19$ &  $[0, 4, 1, 1]~$ & $[4, 4, 1, 1]~$ & $[4, 0, 0, 0]~$ & $[0, 0, 0, 0]~$ & $[0, 0, 1, 1]$ & $[4, 0, 1, 1]$ & $[0, 4, 0, 0]$\\ \hline
$2822 = 2.17.83$ &  $[0, 6, 0]~$ & $[24, 6, 0]$ &$[24, 0, 0]$ & $[0, 6, 0]$ & $[24, 0, 0]$ & $[0, 0, 0]$ & $[24, 6, 0]$\\ \hline
$5654 = 2.257.11$ &  $[0, 4, 0]~$ & $[112, 0, 0]~$ & $[112, 4, 0]~$ & $[112, 0, 0]~$ & $[0, 0, 0]$ & $[112, 4, 0]~$ & $[0, 4, 0]$\\ \hline
$8854 = 2.233.19$  & $[60, 4, 1, 1]~$ & $[0, 4, 0, 0]~$ & $[60, 0, 1, 1]~$ & $[60, 0, 1, 1]~$ & $[60, 4, 1, 1]~$ & $[0, 4, 0, 0]~$ & $[0, 0, 0, 0]~$\\ \hline
$10806 = 2.1801.3$ & $[24, 24, 0,1]~$ & $[24, 24, 0, 0]~$ & $[0, 0, 0, 1]~$ & $[0, 0, 0, 0]~$ & $[24, 24, 0, 0]~$ & $[24, 24, 0, 1]~$ & $[0, 0, 0, 1]~$\\ \hline
\end{longtable}
}
(ii) Assume $p(x\pm1)$ is a square in $\NN$, hence, according to the Proposition $\ref{28}$, we have  $E_{\KK_2}=\langle i, \sqrt{i\epsilon_{q}}, \sqrt{i\epsilon_{2p}}, \sqrt{\epsilon_{2pq}}\rangle$. Thus  proceeding as in the case (i) we prove that $\h$,  $\o$ and $\o\h$ do not capitulate in  $\KK_2$. The following examples illustrate these results.\\
\indent (a) First case: $p(x+1)$ is a square in $\NN$.
The first table gives examples of the ideals $\mathcal{I}$,  $\o$ and  $\h$ which are not principal in $\kk$, and gives the structures of the class groups of $\kk$ and $\KK_2$ respectively; whereas the second table gives the cases of  capitulation  of these ideals in $\KK_2$.
{\tiny
\begin{longtable}{| p{0.77in} | p{0.45in} | p{0.35in}| p{0.35in} |p{0.4in} |p{0.4in} | p{0.4in} | p{0.6in} |}
\hline
$d = 2pq$  & $p(x+1)$ & $\mathcal{I}$& $\mathcal{I}^2$ & $\o$ & $\h$ &  $\mathbf{C}l(\kk)$ & $\mathbf{C}l(\KK_2)$\\
\hline
\endfirsthead
\hline
 $d = 2pq$  & $p(x+1)$ & $\mathcal{I}$ & $\mathcal{I}^2$ & $\o$ & $\h$  & $\mathbf{C}l(\kk)$ & $\mathbf{C}l(\KK_2)$\\
\hline
\endhead
$3358 = 2.73.23$ & $217248^2$ & $[4, 0, 0]~$ & $[0, 0, 0]~$ & $[0, 2, 1]~$ & $[0, 2, 0]~$ & $(8, 4, 2)$ & $(96, 8, 2, 2)$\\ \hline
$3502 = 2.17.103$ & $ 447916^2$ & $[2, 2, 0]~$ & $[0, 0, 0]~$ & $[2, 0, 1]~$ & $[0, 2, 0]~$ & $(4, 4, 2)$ & $(20, 4, 2, 2, 2)$ \\ \hline
$6014 = 2.97.31$ & $388^2$ & $[12, 0, 0]$ & $[0, 0, 0]~$ & $[0, 4, 1]$ & $[12, 4, 0]~$ & $(24, 8, 2)$ & $(240, 24, 2, 2)$\\  \hline
$9118 = 2.97.47$ & $11181384^2$ & $[4, 0, 0]~$ & $[0, 0, 0]~$ & $[4, 0, 1]~$ & $[0, 2, 0]~$ & $(8, 4, 2)$ & $(20, 20, 4, 2, 2)$\\ \hline
\end{longtable}
\begin{longtable}{| p{0.4in}  | p{0.54in} |p{0.55in} |p{0.45in} | p{0.5in} | p{0.55in} |p{0.54in} |p{0.45in} |}
\hline
 $d = 2pq$  & $\o\mathcal{O}_{\KK_2}$ & $\h\mathcal{O}_{\KK_2}$ & $\o\h\mathcal{O}_{\KK_2}$ & $\mathcal{I}\mathcal{O}_{\KK_2}$ & $\h\mathcal{I}\mathcal{O}_{\KK_2}$ & $\o\mathcal{I}\mathcal{O}_{\KK_2}$ & $\o\h\mathcal{I}\mathcal{O}_{\KK_2}$\\
\hline
\endfirsthead
\hline
 $d = 2pq$  & $\o\mathcal{O}_{\KK_2}$ & $\h\mathcal{O}_{\KK_2}$ & $\o\h\mathcal{O}_{\KK_2}$ & $\mathcal{I}\mathcal{O}_{\KK_2}$ & $\h\mathcal{I}\mathcal{O}_{\KK_2}$ & $\o\mathcal{I}\mathcal{O}_{\KK_2}$ & $\o\h\mathcal{I}\mathcal{O}_{\KK_2}$\\
\hline
\endhead
$3358 = 2.73.23$ & $[48, 4, 0, 0]~$ & $[48, 0, 0, 0]~$ & $[0, 4, 0, 0]$ & $[0, 4, 0, 0]$ & $[48, 4, 0, 0]$ & $[48, 0, 0, 0]$ & $[0, 0, 0, 0]$\\ \hline
$3502 = 2.17.103$ & $[0, 0, 1, 0, 0]~$ & $[0, 2, 0, 0, 0]~$  & $[0, 2, 1, 0, 0]~$ &  $[0, 2, 0, 0, 0]~$ & $[0, 0, 0, 0, 0]~$ & $[0, 2, 1, 0, 0]$ & $[0, 0, 1, 0, 0]$\\ \hline
$6014 = 2.97.31$ & $[120, 12, 0, 0]$ & $[120, 0, 0, 0]~$ & $[0, 12, 0, 0]$ & $[120, 12, 0, 0]~$ & $[0, 12, 0, 0]$ & $[0, 0, 0, 0]$ & $[120, 0, 0, 0]$\\ \hline
$9118 = 2.97.47$ & $[10, 10, 2, 1, 0]~$ & $[10, 10, 2, 0, 0]~$ & $[0, 0, 0, 1, 0]~$ & $[0, 0, 0, 0, 0]$ & $[10, 10, 2, 0, 0]~$ & $[10, 10, 2, 1, 0]~$ & $[0, 0, 0, 1, 0]$ \\ \hline
\end{longtable}
}
 (b) Second  case: $p(x-1)$ is a square in $\NN$.
The first table gives examples of  the ideals $\mathcal{I}$,  $\o$ and  $\h$ which are not principal in $\kk$, and gives the structures of the class groups of $\kk$ and $\KK_2$ respectively; whereas the second table gives the cases of  capitulation  of these ideals in $\KK_2$.
{\tiny
\begin{longtable}{| p{0.84in} | p{0.42in} | p{0.4in} | p{0.4in}|p{0.4in} |p{0.4in} | p{0.4in} | p{0.5in} |}
\hline
$d = 2pq$  & $p(x-1)$ & $\mathcal{I}$ & $\mathcal{I}^2$ & $\o$ & $\h$ &  $\mathbf{C}l(\kk)$ & $\mathbf{C}l(\KK_2)$\\
\hline
\endfirsthead
\hline
 $d = 2pq$  & $p(x-1)$ & $\mathcal{I}$ & $\mathcal{I}^2$ & $\o$ & $\h$  & $\mathbf{C}l(\kk)$ & $\mathbf{C}l(\KK_2)$\\
\hline
\endhead
$438 = 2.73.3$ & $ 21316$ & $[0, 1, 1]~$ & $[0, 0, 0]~$ & $[2, 1, 1]~$ & $[2, 0, 0]~$ & $(4, 2, 2)$ & $(32, 2, 2, 2)$\\ \hline
$2022 = 2.337.3$ & $454276$ & $[6, 1, 0]~$  & $[0, 0, 0]~$ & $[0, 0, 1]~$ & $[6, 0, 0]~$ & $(12, 2, 2)$ & $(48, 24, 2)$ \\ \hline
$2598 = 2.433.3$ & $749956$ & $[6, 1, 1]~$ & $[0, 0, 0]~$ & $[0, 1, 1]~$ & $[6, 0, 0]~$ & $(12, 2, 2)$ & $(132, 4, 4)$\\  \hline
$5622 = 2.937.3$ & $3511876$ & $[0, 2, 1]~$ & $[0, 0, 0]~$ & $[0, 0, 1]~$ & $[8, 2, 0]~$ & $(16, 4, 2)$ & $(224, 8, 4)$\\ \hline
\end{longtable}
\begin{longtable}{| p{0.711in}  | p{0.43in} |p{0.43in} |p{0.45in} | p{0.42in} | p{0.4in} |p{0.4in} |p{0.45in} |}
\hline
 $d = 2pq$  & $\o\mathcal{O}_{\KK_2}$ & $\h\mathcal{O}_{\KK_2}$ & $\o\h\mathcal{O}_{\KK_2}$ & $\mathcal{I}\mathcal{O}_{\KK_2}$ & $\h\mathcal{I}\mathcal{O}_{\KK_2}$ & $\o\mathcal{I}\mathcal{O}_{\KK_2}$ & $\o\h\mathcal{I}\mathcal{O}_{\KK_2}$\\
\hline
\endfirsthead
\hline
 $d = 2pq$  & $\o\mathcal{O}_{\KK_2}$ & $\h\mathcal{O}_{\KK_2}$ & $\o\h\mathcal{O}_{\KK_2}$ & $\mathcal{I}\mathcal{O}_{\KK_2}$ & $\h\mathcal{I}\mathcal{O}_{\KK_2}$ & $\o\mathcal{I}\mathcal{O}_{\KK_2}$ & $\o\h\mathcal{I}\mathcal{O}_{\KK_2}$\\
\hline
\endhead
$438 = 2.73.3$ & $[16, 1, 1, 1]~$ & $[0, 1, 1, 1]~$ & $[16, 0, 0, 0]~$ & $[0, 1, 1, 1]~$ & $[0, 0, 0, 0]$ & $[16, 0, 0, 0]~$ & $[16, 1, 1, 1]~$\\ \hline
$2022 = 2.337.3$ & $[24, 12, 0]~$ & $[0, 12, 0]~$ & $[24, 0, 0]~$ & $[24, 12, 0]~$ & $[24, 0, 0]~$ & $[0, 0, 0]~$ & $[0, 12, 0]~$\\ \hline
$2598 = 2.433.3$ & $[66, 2, 0]~$ & $[0, 2, 2]~$ & $[66, 0, 2]~$ & $[66, 0, 2]~$ & $[66, 2, 0]$ & $[0, 2, 2]$ & $[0, 0, 0]$\\ \hline
$5622 = 2.937.3$ & $[112, 4, 0]~$ & $[112, 0, 0]~$ & $[0, 4, 0]~$ & $[0, 0, 0]~$ & $[112, 0, 0]~$ & $[112, 4, 0]~$ & $[0, 4, 0]~$ \\ \hline
\end{longtable}
}
 \indent (3) Suppose that $N(\epsilon_{2p})=-1$.
 Prove that  $\o\h$ and $\o\hh$ capitulate in $\KK_2$.
Put  $\epsilon_{2p}=a+b\sqrt{2p}$, then $a^2+1=2b^2p$, hence by the decomposition uniqueness in  $\ZZ[i]$ there exist  $b_1$ and $b_2$ in $\ZZ[i]$ such that:
 $$ \left\{\begin{array}{ll}
a\pm i &=b_1^2(1+i)\pi_1,\\
a\mp i &= b_2^2(1-i)\pi_2,
\end{array}\right. \text{ or }
  \left\{\begin{array}{ll}
a\pm i &=i(1+i)b_1^2\pi_1,\\
a\mp i &=-i(1-i)b_2^2\pi_2,
\end{array}\right.\text{ with  }b=b_1b_2.$$
 Consequently $\sqrt{\epsilon_{2p}}=\frac{1}{2}(b_1(1+i)\sqrt{(1\pm i)\pi_1}+b_2(1-i)\sqrt{(1\mp i)\pi_2})$, hence $(1\pm i)\pi_1\epsilon_{2p}$ and $(1\mp i)\pi_2\epsilon_{2p}$ are squares in $\KK_2$. Thus
 $(\alpha^2)= ((1\pm i)\pi_1)$ and $(\beta^2)=((1\mp i)\pi_2)$, with some $\alpha$, $\beta$  in $\KK_2$. Therefore $\o\h=(\alpha)$ and $\o\hh=(\beta)$ i.e.  $\o\h$ and $\o\hh$ capitulate in $\KK_2$.\\
\indent (i) If $x\pm1$ is a square in $\NN$, then Lemma \ref{7} yields that   $\h\hh$,   $\o\h$ and $\o\hh$ are not principal in $\kk$, hence the result.\\
\indent (ii) If $x+$ and $x-1$  are not squares in  $\NN$, then Lemma \ref{7} yields that $[\o\h]=[\o\hh]$, hence the result.
\end{proof}
\subsection{Capitulation in $\KK_3$} Let $\KK_3=\kk(\sqrt 2)=\QQ(\sqrt 2, \sqrt{pq}, i)$ and put $\epsilon_{pq}=a+b\sqrt{pq}$, $\epsilon_{2pq}=x+y\sqrt{2pq}$. Let $Q_{\KK_3}$ denote the unit index of  $\KK_3$.
\begin{The}\label{231}
Keep  the notations and hypotheses previously  mentioned.
\begin{enumerate}[\rm\indent(1)]
\item If  both of $x\pm1$ and  $a\pm1$  are squares in $\NN$, then
\begin{enumerate}[\rm(a)]
\item If $Q_{\KK_3}=2$, then $\kappa_{\KK_3}=\langle[\o]\rangle$.
\item If  $Q_{\KK_3}=1$, then $\kappa_{\KK_3}=\langle[\o], [\h\hh]\rangle$.
\end{enumerate}
\item If   $x\pm1$ is a square in $\NN$ and $a+1$, $a-1$  are not, then $\kappa_{\KK_3}=\langle[\o], [\h\hh]\rangle$.
\item If  $a\pm1$ is a square in $\NN$ and  $x+1$, $x-1$  are not, then there exist an unambiguous ideal $\mathcal{I}$ in  $\kk/\QQ(i)$ of order $2$ such that $\kappa_{\KK_3}=\langle[\o], [\mathcal{I}]\rangle$ or
$\langle[\o],  [\h\mathcal{I}]\rangle$.
\item If $x+1$, $x-1$, $a+1$ and $a-1$ are note squares in $\NN$, then $\kappa_{\KK_3}=\langle[\o]\rangle$.
\end{enumerate}
\end{The}
\begin{proof}
As $N(\epsilon_2)=-1$, then  $\sqrt{(1+i)\epsilon_2}=\frac{1}{2}(2+(1+i)\sqrt2)$. Hence there exist  $\beta\in\KK_3$ such that $\mathcal{H}_0^2=(1+i)=(\beta^2)$, therefore $\mathcal{H}_0$ capitulates in $\KK_3$.\\
\indent (1) Assume $x\pm1$ and  $a\pm1$ are squares in $\NN$.\\
 (a) If  $Q_{\KK_3}=2$, then by the Theorem \ref{226},  $|\kappa_{\KK_3}|=2$, hence $\kappa_{\KK_3}=\langle[\o]\rangle$.\\
 (b) If  $Q_{\KK_3}=1$, then  by the Theorem \ref{226},  $|\kappa_{\KK_3}|=4$. Since   $a\pm1$ is a square in $\NN$, so the Lemma \ref{247} yields that $p\equiv1\pmod8$. Therefore the Proposition \ref{248} implies that
  $$\mathrm{Am}(\kk/\QQ(i))=\mathrm{Am}_s(\kk/\QQ(i))=\langle[\o], [\h], [\hh]\rangle.$$
  Proceeding as in the proof of Theorem \ref{229}(\ref{item2}), we show that $\h$ and  $\hh$ do not capitulate in $\KK_3$. On the other hand, as
    $|\kappa_{\KK_3}|=4$ and $\kappa_{\KK_3}\subseteq \mathrm{Am}(\kk/\QQ(i))$, so necessarily   $\h\hh$ capitulate in $\KK_3$. Finally, Lemma \ref{7} yields that  $\h\hh$, $\o$ and  $\o\h\hh$ are not principal in  $\kk$. Thus the result.\\
\indent (2)  Assume    $x\pm1$ is a square in $\NN$ and $a+1$, $a-1$  are not. As $\o$ capitulates in $\KK_3$ and   $|\kappa_{\KK_3}|=4$ (Theorem \ref{226}), it suffices to prove that $\h\hh$ capitulates in $\KK_3$. According to  the proof of Proposition \ref{29}, we have $p\epsilon_{pq}$  is a square in $\KK_3$; hence there exist
  $\alpha$  in $\KK_3$ such that $(p)=(\alpha^2)$, so $\h\hh=(\alpha)$. Thus the result.\\
\indent (3) If  $a\pm1$ is a square in $\NN$ and $x+1$, $x-1$  are not, then Lemma \ref{247} implies that $p\equiv1\pmod8$; hence we are in the hypotheses of the Proposition \ref{248}. On the other hand, from the Lemma \ref{7} we get  $[\h]=[\hh]$. Therefore, proceeding as in the proof of Theorem  \ref{229}, we show that  $\h$ does not capitulate in $\KK_3$. The following examples clarify the two cases of capitulation:
{\footnotesize
\begin{longtable}{| c | c | c | c | c | c | c |}
\hline
  $d = 2pq$ & $\o\mathcal{O}_{\KK_3}$ & $\h\mathcal{O}_{\KK_3}$ & $\mathcal{I}\mathcal{O}_{\KK_3}$ & $\h\mathcal{I}\mathcal{O}_{\KK_3}$ &  $\mathbf{C}l(\kk)$ & $\mathbf{C}l(\KK_2)$\\
\hline
\endfirsthead
\hline
 $d = 2pq$ & $\o\mathcal{O}_{\KK_3}$ & $\h\mathcal{O}_{\KK_3}$ & $\mathcal{I}\mathcal{O}_{\KK_3}$ & $\h\mathcal{I}\mathcal{O}_{\KK_3}$ &  $\mathbf{C}l(\kk)$ & $\mathbf{C}l(\KK_2)$\\
\hline
\endhead
 $582 = 2.97.3$ & $[0, 0, 0]$ & $[8, 2, 0]$ &  $[0, 0, 0]$  & $[8, 2, 0]$ & $(8, 2, 2)$ & $(16, 4, 2)$\\ \hline
 $2006 = 2.17.59$ & $[0, 0, 0]$ & $[24, 0, 0]$ &  $[24, 0, 0]$  & $[0, 0, 0]$ & $(24, 2, 2)$ & $(48, 4, 2)$\\ \hline
 $2454 = 2.409.3$ & $[0, 0, 0]$ & $[16, 0, 0]$ & $[16, 0, 0]$ & $[0, 0, 0]$ & $(16, 2, 2)$ & $(32, 4, 2)$\\ \hline
 $2742 = 2.457.3$ & $[0, 0, 0]~$ & $[48, 2, 0]~$ & $[0, 0, 0]~$ & $[48, 2, 0]$ & $ (16, 2, 2) $ & $ (96, 4, 2) $\\ \hline
\end{longtable}
}
\indent (4) Suppose that $x+1$, $x-1$,   $a+1$ and $a-1$ are not squares in $\NN$, then $|\kappa_{\KK_3}|=2$ (Theorem \ref{226}). Thus
$\kappa_{\KK_3}=\langle[\o]\rangle$.
\end{proof}
From theorems \ref{227}, \ref{229} and \ref{231}  we deduce the following theorem.
\begin{The}\label{9}
Let $\kk=\QQ(\sqrt{2pq},i)$, where $p\equiv-q\equiv1 \pmod 4$ are different primes,  and $\G$ its genus field. Put   $\epsilon_{2pq}=x+y\sqrt{2pq}$ and $\epsilon_{pq}=a+b\sqrt{pq}$.
\begin{enumerate}[\rm\indent(1)]
  \item If $x\pm1$ is a square in $\NN$, then $\langle[\o], [\h], [\hh]\rangle\subseteq  \kappa_{\k}$.
  \item If $x+1$ and $x-1$ are not squares in $\NN$, then
\begin{enumerate}[\rm(a)]
\item If  $N(\epsilon_{2p})=1$ or $a\pm1$ is a square in  $\NN$, then there exists an unambiguous ideal  $\mathcal{I}$ in $\kk/\QQ(i)$ of order $2$ such that: $\langle[\o], [\h], [\mathcal{I}] \rangle\subseteq \kappa_{\k}.$
\item Else
  $\langle[\o], [\h] \rangle\subseteq \kappa_{\k}$.
\end{enumerate}
\end{enumerate}
\end{The}
 Theorem \ref{9} implies the following corollary:
\begin{coro}
Let $\kk=\QQ(\sqrt{2pq},i)$, where $p\equiv-q\equiv1 \pmod 4$ are different primes.   Let $\G$ be the genus field of $\kk$ and  $\mathrm{A}m_s(\kk/\QQ(i))$ be the group of  the strongly ambiguous class  of $\kk/\QQ(i)$, then $\mathrm{A}m_s(\kk/\QQ(i))\subseteq \kappa_{\k}$.
\end{coro}
\section{\bf{Application}}
Let $p\equiv -q\equiv1\pmod4$ be different primes such that $p\equiv1\pmod8$,  $q\equiv3\pmod8$ and $\left(\frac{p}{q}\right)=-1$. Hence, according to  \cite{AzTa-08},  $\mathbf{C}l_2(\mathds{\kk})$ is of type $(2 ,2, 2)$.   Therefore, under these assumptions,    $\mathbf{C}l_2(\mathds{\kk})=\mathrm{A}m_s(\kk/\QQ(i))=\langle [\mathcal{H}_0], [\mathcal{H}_1], [\mathcal{H}_2]\rangle$ (see \cite{AZT12-2}). To continue we need the following result.
\begin{lem}\label{8}
Let $\epsilon_{2pq}=x+y\sqrt{2pq}$ $($resp. $\epsilon_{pq}=a+b\sqrt{pq})$ denote the fundamental unit of $\QQ(\sqrt{2pq})$ $($resp. $\QQ(\sqrt{pq}))$. Then
\begin{enumerate}[\rm\indent(1)]
  \item  $x-1$ is a square in $\NN$.
  \item $a-1$ is a square in $\NN$.
\end{enumerate}
\end{lem}
\begin{proof} (1)  By Lemma \ref{5} and according to the decomposition uniqueness in $\ZZ$, there are  six cases to discus: $x\pm1$ or $p(x\pm1)$ or $2p(x\pm1)$ is a square in $\NN$.
\begin{enumerate}[a.]
\item If $x+1$ is a square in $\NN$, then $\left\{\begin{array}{ll}
  x+1=y_1^2,\\
  x-1=2pqy_2^2,
  \end{array}\right.$\\ hence $1=\left(\frac{x+1}{q}\right)=\left(\frac{x-1+2}{q}\right)=\left(\frac{2}{q}\right)$, which contradicts  the fact that $\left(\frac{2}{q}\right)=-1$.
\item If $p(x\pm1)$ is a square in $\NN$, then
  $\left\{\begin{array}{ll}
  x\pm1=py_1^2,\\
  x\mp1=2qy_2^2,
  \end{array}\right.$\\  hence $\left(\frac{2q}{p}\right)=\left(\frac{x\mp1}{p}\right)=\left(\frac{x\pm1\mp2}{p}\right)=\left(\frac{2}{p}\right)$, thus $\left(\frac{q}{p}\right)=1$. Which is false, since  $\left(\frac{p}{q}\right)=-1$.
\item If $2p(x+1)$ is a square in $\NN$, then
  $\left\{\begin{array}{ll}
  x+1=py_1^2,\\
  x-1=2qy_2^2,
  \end{array}\right.$\\ hence $\left(\frac{2p}{q}\right)=\left(\frac{x+1}{q}\right)=\left(\frac{x-1+2}{q}\right)=\left(\frac{2}{q}\right)$, which leads to the contradiction  $\left(\frac{q}{p}\right)=1$.
\item If $2p(x-1)$ is a square in $\NN$, then
  $\left\{\begin{array}{ll}
  x-1=py_1^2,\\
  x+1=2qy_2^2,
  \end{array}\right.$\\ hence $\left(\frac{q}{p}\right)=\left(\frac{x+1}{p}\right)=\left(\frac{x-1+2}{p}\right)=\left(\frac{2}{p}\right)=1$, which is false.
\end{enumerate}
 Consequently, the only case which is possible is: $x-1$ is a square in $\NN$.

 \indent (2) Proceeding similarly, we show that   $a-1$ is a square in $\NN$.
\end{proof}

\begin{The}
Let  $\kk=\QQ(\sqrt{2pq}, i)$, where $p\equiv -q\equiv1\pmod4$ are different primes satisfying the conditions $p\equiv1\pmod8$,  $q\equiv3\pmod8$ and $\left(\frac{p}{q}\right)=-1$. Put
 $\KK_1=\kk(\sqrt {p})$, $\KK_2=\kk(\sqrt {q})$ and $\KK_3=\kk(\sqrt {2})$. Let $\k$ denote the absolute genus field of $\kk$ and $(\kk/\QQ(i))^*$ its  relative genus field over $\QQ(i)$.
\begin{enumerate}[\rm\indent(1)]
\item $\k\varsubsetneq (\kk/\QQ(i))^*$.
  \item $\kappa_{\KK_1}=\langle[\mathcal{H}_1], [\mathcal{H}_2]\rangle$.
  \item  Denote by $\epsilon_{2p}$  the fundamental unit of $\QQ(\sqrt{2p})$, so:
\begin{enumerate}[\rm(a)]
  \item If  $N(\epsilon_{2p})=1$, then $\kappa_{\KK_2}=\langle[\mathcal{H}_1], [\mathcal{H}_2]\rangle$ or $\langle[\o], [\h\hh]\rangle$.
  \item Else,  $\kappa_{\KK_2}=\langle[\o\h], [\o\hh]\rangle$
\end{enumerate}
  \item Denote by $Q_{\KK_3}$ the unit index of  $\KK_3$, so:
 \begin{enumerate}[\rm(a)]
  \item If $Q_{\KK_3}=1$, then $\kappa_{\KK_3}=\langle[\o], [\h\hh]\rangle$.
  \item If $Q_{\KK_3}=2$, then $\kappa_{\KK_3}=\langle[\o]\rangle$
 \end{enumerate}
  \item $\kappa_{\G}=\mathrm{A}m_s(\kk/\QQ(i))=\mathrm{C}l_2(\kk)$.
\end{enumerate}
\end{The}
\begin{proof}
\begin{enumerate}[\noindent(1)]
\item From  Lemma \ref{8}, we have $x-1$ is a square in $\NN$. Then  Proposition \ref{248} yields the first assertion.
\item From  Lemma \ref{8}, we have $x-1$ is a square in $\NN$. Then  Theorem \ref{227}(1)  yields the second assertion.
\item From  Lemma \ref{8}, we have $x-1$ is a square in $\NN$. Therefore
\begin{enumerate}
  \item If  $N(\epsilon_{2p})=1$, then Theorem \ref{229}(1) yields the result.
  \item If  $N(\epsilon_{2p})=-1$, then Theorem \ref{229}(3) yields the result.
\end{enumerate}
\item As $x-1$ and $a-1$ are  squares in $\NN$ (Lemma \ref{8}), so  Theorem \ref{231}(1) yields the result.
\item As $p\equiv1\pmod8$, so from Proposition \ref{248} we get $\mathrm{A}m_s(\kk/\QQ(i))=\langle[\o], [\h], [\hh]\rangle$. Hence $\mathrm{A}m_s(\kk/\QQ(i))=\mathrm{C}l_2(\kk)$. The assertions (2), (3) and (4) imply that $\kappa_{\G}=\mathrm{A}m_s(\kk/\QQ(i))=\mathrm{C}l_2(\kk)$.
\end{enumerate}
\end{proof}
\section*{Acknowledgement}
We would like to thank the unknown referee  for his several helpful suggestions and for calling our attention
to the missing details.

{\small
{\em Authors' addresses}:\\
{\em Abdelmalek Azizi and Abdelkader Zekhnini}, Mohammed First University, Mathematics Department, Sciences Faculty, Oujda, Morocco.\\
 e-mail: \texttt{abdelmalekazizi@\allowbreak yahoo.fr}\\
e-mail: \texttt{zekha1@\allowbreak yahoo.fr}.\\
{\em Mohammed Taous}, Moulay Ismail University, Mathematics Department, Sciences and Techniques Faculty, Errachidia, Morocco.\\
 e-mail: \texttt{taousm@\allowbreak hotmail.com}
}

\end{document}